\documentclass[11pt,reqno]{amsart}
\usepackage{newlfont,amsmath,amsthm,amsgen,amscd}
\usepackage{amssymb}
\usepackage{latexsym}
\usepackage{fancyhdr}
\usepackage{ifthen}
\usepackage{amsfonts}
\usepackage{lastpage}
\usepackage{afterpage}
\usepackage{epsfig}
\usepackage{enumerate}
\usepackage{euscript,color}

\newcommand{\belabel}[1]{\begin{equation}\label{#1}}
\newtheorem{mthm}{Theorem}
\newtheorem{mcor}{Corollary}

\newcommand{\ph}{\varphi}

\newcommand{\bM}{{\overline{\cal M}}}

\newcommand{\bg}{{\overline{\mathbf{g}}}}

\newcommand{\fX}{{\mathfrak{X}}}

\newcommand{\Id}{{\mathrm{Id}}}
\newcommand{\F}{{\cal F}}
\newcommand{\beq}{\begin{eqnarray*}}
\newcommand{\eeq}{\end{eqnarray*}}
\newcommand{\be}{\begin{eqnarray}}
\newcommand{\ee}{\end{eqnarray}}

\newcommand{\beqn}{\begin{equation}}
\newcommand{\eeqn}{\end{equation}}

\theoremstyle{definition}

\newtheorem{re}{Remark}[section]

\newtheorem{bsp}{Example}[section]

\newtheorem*{bsp*}{Example}
\newtheorem*{def*}{Definition}
\theoremstyle{plain}
\newtheorem{Lemma}{Lemma}[section]
\newtheorem*{lem*}{Lemma}
\newtheorem{Proposition}{Proposition}[section]
\newtheorem{Corollary}{Corollary}[section]
\newtheorem{Theorem}{Theorem}[section]
\newtheorem*{theo*}{Theorem}
\newtheorem*{conj*}{Conjecture}

\newcommand{\M}{{\cal M}}
\newcommand{\cF}{{\cal F}}
\newcommand{\g}{\mathbf{g}}

\newcommand{\W}{\mathrm{W}}

\newcommand{\del}{\partial}

\newcommand{\bleml}[1]{\begin{Lemma} \label{#1}}
\newcommand{\blem}{\begin{Lemma}}
\newcommand{\elem}{\end{Lemma}}

\newcommand{\btheo}{\begin{Theorem}}
\newcommand{\btheol}[1]{\begin{Theorem}\label{#1}}
\newcommand{\etheo}{\end{Theorem}}

\newcommand{\bpropl}[1]{\begin{Proposition} \label{#1}}
\newcommand{\bprop}{\begin{Proposition}}
\newcommand{\eprop}{\end{Proposition}}

\newcommand{\bcorl}[1]{\begin{Corollary} \label{#1}}
\newcommand{\bcor}{\begin{Corollary}}
\newcommand{\ecor}{\end{Corollary}}

\newcommand{\bbem}{\begin{re}}
\newcommand{\ebem}{\end{re}}

\newcommand{\bprf}{\begin{proof}}
\newcommand{\eprf}{\end{proof}}

\begin{document}
\title[The Cauchy problem for parallel spinors as $1^{\text{st}}$-order symmetric hyperbolic system]{The Cauchy problem for parallel spinors as $1^{\text{st}}$-order symmetric hyperbolic system}


\author{Andree Lischewski} 
\address[Andree Lischewski]{School of Mathematical Sciences, University of Adelaide, SA~5005,
Australia} \email{andree.lischewski@adelaide.edu.au}
\thanks{The
 author acknowledges support from the Collaborative Research Center
647 ``Space-Time-Matter'' of the German Research Foundation.}
\begin{abstract}
We prove that a smooth Riemannian manifold admitting an imaginary generalized Killing spinor whose Dirac current satisfies an additional algebraic constraint condition can be embedded as spacelike Cauchy hypersurface in a smooth Lorentzian manifold on which the given spinor extends to a null parallel spinor. 
%This generalizes results obtained before for the analytic category and 
This is in contrast to a corresponding Cauchy problem for real generalized Killing spinors into Riemannian manifolds.
The construction is based on first order symmetric hyperbolic PDE-methods. In fact, the coupled evolution equations for metric and spinor as considered here extend and generalize the well known PDE-system appearing in the Cauchy problem for the vacuum Einstein equations.
Special cases are discussed and the statement is compared with a similar result obtained recently for the analytic category.
%On a  Lorentzian manifold the existence of a parallel null vector field implies certain constraint conditions on the induced Riemannian geometry of a space-like hypersurface. We will derive these constraint conditions and, conversely, show that every real analytic Riemannian manifold satisfying the constraint conditions can be extended to a Lorentzian manifold  with a parallel null vector field. Similarly, every parallel null spinor on a Lorentzian manifold induces an imaginary generalized Killing spinor on a space-like hypersurface. Then, based on  the fact that a parallel  spinor field induces a parallel  vector field, we can apply the first result to prove: every real analytic Riemannian manifold carrying a real analytic, imaginary generalized Killing spinor can be extended to a Lorentzian manifold with a parallel null spinor. Finally, we give examples of geodesically complete Riemannian manifolds satisfying the constraint conditions.
\end{abstract}

\maketitle

%\tableofcontents

\section{Background and  main result}

%\subsection*{Special Lorentzian holonomy and the related Cauchy problems}
%This paper is a contribution to the research programme of studying global and causal properties of Lorentzian manifolds with special holonomy. A Lorentzian manifold has {\em special holonomy} if  the connected component of its holonomy group is reduced from the full group $\mathbf{SO}^0(1,n)$, but still acts indecomposably, i.e., without non-degenerate invariant subspaces. In this situation the Lorentzian manifold admits a  bundle of tangent null lines that is invariant  under parallel transport. The possible special Lorentzian holonomy groups were classified in \cite{bb-ike93} and \cite{leistnerjdg}, all of them can be realised by local metrics \cite{galaev05}, but many questions about the consequences of special holonomy for global and causal properties of the manifold are still open.  A special case of this situation is when the parallel null line bundle is spanned by a {\em parallel null vector field}. This is the case we will study in this paper. It  is motivated by the question which Lorentzian manifolds admit a {\em parallel spinor field}, which in turn draws its motivation from mathematical physics.

Let $(M,g)$ be a smooth, time-oriented Lorentzian spin manifold with complex spinor bundle $S=S^g \rightarrow M$ and spinor covariant derivative $\nabla = \nabla^{S}$ induced by the Levi Civita connection of $g$. Moreover, we assume that there exists a nontrivial parallel spinor on $(M,g)$, i.e. a solution $\phi \in \Gamma(M,S^g)$ of $\nabla \phi = 0$. This overdetermined PDE and its relation to Lorentzian geometries naturally arises in various areas of mathematical physics, for instance in supergravity since the 1970s, cf. \cite{Deser,sugr}. In fact, in various dimensions the classification and construction of Lorentzian manifolds admitting parallel spinors is directly linked to the construction of supersymmetric supergravity backgrounds (with zero flux), as for instance demonstrated for M-theory backgrounds in dimension 11 in \cite{brmwa}. One construction principle for supersymmetric M-theory backgrounds with zero flux is motivated by the Cauchy problem for the Einstein equations in classical general relativity (cf. \cite{choquet,rin} for an overview). 
Let us to this end additionally assume that $\Sigma \subset M$ is a smooth spacelike hypersurface equipped with the induced metric $g^{\Sigma}$ and spin structure and Weingarten tensor or second fundamental form $W$. It is well known that the vacuum Einstein field equation $\text{Ric}^g = 0$ imposes on $(\Sigma,g^{\Sigma})$ the \textit{constraint equations}
\begin{equation}
\begin{aligned} \label{con}
d \text{tr}_{g^{\Sigma}}W  + \delta^{\Sigma} W &= 0, \\
\text{scal}^{\Sigma} - \text{tr}_{g^{\Sigma}}(W^2) +  (\text{tr}_{g^{\Sigma}} W)^2 &=0.
\end{aligned}
\end{equation}
Conversely, given $(\Sigma,g^{\Sigma},W)$ solving the system \eqref{con} one can find a Ricci flat and globally hyperbolic Lorentzian metric $(M,g)$ defined on a neighborhood of $\Sigma$ in $\mathbb{R} \times \Sigma$ in which $(\Sigma,g^{\Sigma})$ embeds with $W$, see \cite{rin} and references therein. The metric $g$ is obtained as solution to \textit{evolution equations} which are equivalent to a certain first order symmetric quasilinear hyperbolic system of PDEs as found in \cite{first}.\\
Motivated by this and returning to the supergravity picture, it is natural to study in comparison the Cauchy problem for Lorentzian manifolds admitting parallel spinors. Note that in contrast to the Riemannian case a Lorentzian manifold with parallel spinor is not necessarily Ricci - flat, whence this problem is not a special case of the Cauchy problem for the vacuum Einstein equations. In certain special cases and dimensions, the Cauchy problem in supergravity has already been considered before, cf. \cite{caus}. However, although the construction of supersymmetric supergravity backgrounds with zero flux where 11 is an upper bound for the dimension is a motivation for this work, we allow arbitrary dimension for $M$ and consider the general geometric Cauchy problem for Lorentzian metrics admitting parallel spinors. \\
Concretely, if the data $(M,g,\phi)$ and $(\Sigma,g^{\Sigma},W)$ are chosen as before, it is well-known (cf. \cite{bgm}) that the restriction $\ph$ of $\phi$ to $\Sigma$ satisfies a constraint equation known as the {\em imaginary generalized $W-$Killing spinor equation}
\begin{align} \label{genks}
 \nabla^{\Sigma}_X
\varphi =\frac{i}{2}\, W(X)\cdot \varphi \qquad\forall X\in T\Sigma.
\end{align}
Moreover, we shall assume that the Dirac current $V_{\phi}$ of $\phi$, as to be recalled in section 2, is isotropic. Otherwise, as $V_{\phi}$ is parallel as well and always causal, $(M,g)$ would split into a metric product of a timelike line and a Riemannian factor. From a more physical perspective the condition $g(V_{\phi},V_{\phi})=0$ is in low dimensions needed to ensure that the symmetry superalgebra naturally associated to the supergravity background $(M,g)$ is a Lie superalgebra, cf. \cite{of121,of122}.
It then follows that $V_{\phi} \cdot \phi = 0$. On $\Sigma$, this imposes the additional \textit{algebraic constraint}
\begin{align} \label{algc}
 U_{\varphi} \cdot \varphi = i\,u_{\varphi} \,\varphi, 
\end{align}
where $u_{\ph}=\langle \ph,\ph \rangle$ denotes the length of $\ph$ and $U_{\ph}$ is the Riemannian Dirac current of $\ph$ (as to be made precise in section 2).
In analogy to the vacuum Einstein equations in classical general relativity we prove that also the converse of these observations is true in the following sense:

\begin{mthm} \label{mother}
Let the following data be given:
\begin{enumerate}
\renewcommand{\labelenumi}{(\roman{enumi})}
\item $(\Sigma,g^{\Sigma})$ is a smooth Riemannian spin manifold admitting a nontrivial spinor $\ph$ satisfying the constraint equations \eqref{genks}-\eqref{algc} for some symmetric tensor $W$ on $\Sigma$.
\item $\lambda \in C^{\infty}(\mathbb{R} \times {\Sigma},\mathbb{R}^+)$ is an arbitrary positive smooth function and $h_t$ is any family of Riemannian metrics on $\Sigma$ with $h_0 = g^{\Sigma}$. For any such data we form the Lorentzian metric $h = -\lambda^2 dt^2 + h_t$ (the background metric), defined on an open neighborhood of $\Sigma$ in $\mathbb{R} \times \Sigma$.
\end{enumerate}
Then there exist an open neighborhood $M$ of $\Sigma$ in $\mathbb{R} \times \Sigma$ and a unique smooth Lorentzian metric $g=g_h$ on $M$  such that 
\begin{enumerate}
\item $(M,g)$ is spin and admits a parallel spinor $\phi$ with $g(V_{\phi},V_{\phi}) = 0$, 
\item $\phi$ restricts on $\Sigma$ to $\ph$,
\item  $g_{|\Sigma} = -\lambda^2_{|\Sigma}dt^2 + g^{\Sigma}$,
\item $g_h$ depends on $h$ in terms of the following PDE-system: The contracted difference tensor of the Levi Civita connections of $g$ and $h$ vanishes, i.e.
\begin{align}
E(X) := -\text{tr}_g \left(g(A(\cdot,\cdot),X) \right) = 0 \text{ }\forall X \in TM,
\end{align}
where $A(Y,Z) := \nabla^g_Y Z - \nabla^h_Y Z$ for $Y,Z \in TM$.
\end{enumerate} 
In particular, $(\Sigma,g^{\Sigma})$ embeds into $(M,g)$ with Weingarten tensor $W$. Moreover, for any given $h$, $(M,g_h)$ can be chosen to be globally hyperbolic with spacelike Cauchy hypersurface $\Sigma$, i.e. $\Sigma$ is met by every inextendible timelike curve in $(M,g)$ exactly once.
\end{mthm}
Theorem \ref{mother} shows that the Cauchy problem for smooth Lorentzian manifolds admitting parallel spinors is up to the choice of a background metric well posed. A variety of examples of manifolds solving the constraints \eqref{genks}-\eqref{algc}, including compact ones, can be found in \cite{rad,bam}. \\
\newline
Let us elaborate on the consequences of this result in more detail. 
%As the Dirac current of a parallel spinor $\phi$ is parallel, the holonomy group of the Lorentzian manifold $(M,g)$ constructed via Theorem \ref{mother} fixes a nontrivial lightlike vector.
%Thus, the result is intimately related to the question: What is the intersection of the class of globally
%hyperbolic Lorentzian manifolds with the class of Lorentzian manifolds with special holonomy?  A Lorentzian manifold is {\em globally hyperbolic} if  
%From this point of view, Theorem \ref{mother} opens up a way of constructing globally hyperbolic Lorentzian manifolds with special holonomy if solutions to the constraint equations \eqref{genks}-\eqref{algc} are given. Examples of manifolds solving the constraints can be found in \cite{rad,bam}. In particular, compact examples are known to exist.\\
First, Theorem \ref{mother} is interesting when compared to the analogous Riemannian situation. Namely for Riemannian manifolds, the corresponding Cauchy  problem was studied by Ammann, Moroianu and Moroianu  in \cite{amm} in relation to the Cauchy problem for Ricci-flat manifolds. They show by using the Cauchy-Kowalewski theorem that in the analytic category \textit{real} generalized $W-$Killing spinors can be extended to parallel spinors on Riemannian manifolds. Counterexamples show that this doesn't carry over to the smooth case. It is therefore worth mentioning that Theorem \ref{mother} requires only smoothness of the initial data. This is in analogy to the Cauchy problem for the Einstein equations which can be solved in the smooth case for Lorentzian manifolds (cf. \cite{choquet,rin} for an overview) but only for analytic initial data in the Riemannian situation (cf. \cite{amm,Koiso1981}). \\

%A background metric $h$ such as in Theorem \eqref{mother} also appears for 
The proof of Theorem \ref{mother} is motivated by considering again the Cauchy problem for the vacuum Einstein field equations, which can be tackled by first order hyperbolic methods, cf. \cite{first}.
This idea in mind, we proceed as follows: Given a lightlike parallel spinor $\phi$ on $(M,g)$ it is by using that $\text{Ric}^g$ is isotropic (cf. \cite{ba81}) easy to show that there is a uniquely determined function $f \in C^{\infty}(M)$ such that the data $(g,\phi,f)$ satisfy a coupled PDE of (local) type
\begin{align}
F(\del^2g, \del g, g,  \del \phi , \phi, \del f, f) = 0, \label{Fs}
\end{align}
whose highest order derivatives are  specified by $\text{Ric}^g$, the Ricci tensor of $g$, and $D^g \phi$, the Dirac operator. We observe that for $f= 0$ and $\phi = 0$ the system \eqref{Fs} reduces to the vacuum Einstein equations $\text{Ric}^g = 0$ for which \cite{first} provides first order hyperbolic methods for smooth initial data. Inspired by this, we show that \eqref{Fs} can be rewritten as a symmetric first order quasilinear hyperbolic PDE. A local solution theory for such systems in the smooth setting is available (cf. \cite{tay}). As initial data we choose the generalized $W$-Killing spinor and the hypersurface metric $g^{\Sigma}$. Then the local solutions of \eqref{Fs} produce local data $(g,\phi)$ which are candidates for Theorem \ref{mother}. In the second step of the proof, we show that the properties of $(g,\phi)$ following from $\eqref{Fs}$ already imply $\nabla \phi = 0$. This is achieved by showing that $\nabla_X \phi$ (and other data) lie in the kernel of a suitably constructed normally hyperbolic second order linear differential operator for which the Cauchy problem is (locally) known to be well-posed. Uniqueness of the solution and \eqref{genks}-\eqref{algc} viewed as initial conditions allow to conclude $\nabla \phi = 0$. Then we show that the local data can be patched together to give $(M,g,\phi)$ as in Theorem \ref{mother}.\\
The background metric $h$ from Theorem \ref{mother} enters in \eqref{Fs} (that is, in fact $F=F_h$) in such a way that \eqref{Fs} is indeed equivalent to a first order hyperbolic system. $h$ is used to manipulate the Ricci tensor $\text{Ric}^g$, which is not hyperbolic considered as differential operator acting on the metric. Such a background metric construction, called \textit{hyperbolic reduction}, is also needed for the Cauchy problem in general relativity. Physics literature seems to use exclusively the canonical background metric $h = -dt^2 + g^{\Sigma}$. However, we show that the proof works for a more arbitrary class of such metrics as specified in Theorem \ref{mother}. This later pays off when comparing the Theorem to results obtained in \cite{blli}:
This reference studies the Cauchy problem for Lorentzian manifolds with parallel spinors in the analytic category. If the data $( \Sigma,g^{\Sigma},W,\ph)$ are all assumed to be analytic and if moreover $\lambda$ is any analytic and positive function, the so called  \textit{lapse function}, defined on $\mathbb{R} \times \Sigma$, the Cauchy Kowalewski Theorem can be used to show the existence of a unique family of Riemannian metrics $g_t$ on $\Sigma$ such that $g_{\lambda}=-\lambda^2 dt^2 + g_t$ is an analytic Lorentzian metric on a neighborhood of $\Sigma$ in $\mathbb{R} \times \Sigma$ enjoying the properties (1)-(3) from Theorem \ref{mother}. In the smooth category we do not know whether there is always a background metric $h$ such that the metric $g_h$ is of form $-\widetilde{\lambda}^2 dt^2 + g_t$. We will compare the results from \cite{blli} to Theorem \ref{mother} more carefully in section 4. At this point let us just keep that both approaches produce \textit{families} of Lorentzian metrics with parallel spinors for given initial data, which are parametrized by analytic positive functions respectively background metrics $h$:
\begin{align*}
\lambda : \mathbb{R} \times \Sigma \rightarrow \mathbb{R}^+,  \text{all data analytic} & \stackrel{\text{Ref. [5]}}{\longrightarrow} &  g_{\lambda}=-\lambda^2 dt^2 + g_t\text{ with }\nabla \phi = 0, \\
h \text{ background metric},  \text{all data smooth} & \stackrel{\text{Thm. 1}}{\longrightarrow}  & g_{h} \text{ with }\nabla \phi = 0,  h-\text{dep. via }E^{g,h}=0.\\
\end{align*}

%A proof of Theorem \ref{mother} in the analytic category (without showing that $(M,g)$ is globally hyperbolic, however) can be found in \cite{blli}.\\
%\newline
The article is organized as follows: In section 2 we recall and collect necessary facts about parallel spinors on Lorentzian manifolds and generalized imaginary $W-$Killing spinors on spacelike hypersurfaces. Section 3 is then devoted to the proof of Theorem \ref{mother}. The main ideas and steps for the proof are explained in some length at the beginning. In section 4 we compare our results to those from \cite{blli} and answer under which additional conditions on $(\Sigma, g^{\Sigma},W)$ the development from Theorem \ref{mother} is actually Ricci flat and illustrate this for a nontrivial warped product example.

%\newpage

\section{Relevant facts about parallel- and imaginary generalized Killing spinors}

\noindent Let $(M,g)$ be a time oriented Lorentzian spin manifold of dimension $n+1$ and denote by $(S,\nabla^S)$ its spinor bundle with the covariant derivative induced by the Levi-Civita
connection $\nabla = \nabla^g$. When no confusion is likely to occur we also write $\nabla^S = \nabla$. For a spinor field $\phi$ on $(M,g)$ we define
its Dirac current $V_{\phi} \in \fX(M)$ by
\[ g(V_{\phi},   X) = - \langle  X \cdot \phi, \phi \rangle\,,  \qquad  \forall\ X \in \fX(M).\]
The vector field $V_{\phi}$ is future-oriented, causal, i.e., 
$g(V_{\phi},V_{\phi}) \leq 0$ and the zero sets of $V_{\phi}$ and
$\phi$ coincide. If $\phi$ is parallel, $V_{\phi}$ is parallel as
well, and thus either null or timelike. We call a spinor field $\phi$
{\em null}, if its Dirac current $V_{\phi}$ is null. In this case we
have $V_{\phi} \cdot \phi=0$ and $\langle \phi,\phi \rangle = 0$.
We also introduce the spin Dirac operator 
\[ D = D^g \stackrel{\text{loc.}}{=} \sum_{a=0}^n \epsilon_a s_a \cdot \nabla_{s_a} \] 
of $(M,g)$, where $(s_0,...,s_n)$ denotes a local pseudo-orthonormal basis for $(M,g)$ and $\epsilon_a = g(s_a,s_a)$.

\begin{Lemma} \label{eve}
Let $\phi \in \Gamma(M,{S})$ be a null parallel spinor on a time-oriented Lorentzian manifold $(M,g)$ with Dirac current $V= V_{\phi}$. Then there exists a function $f \in C^{\infty}(M)$ such that\footnote{This is the system mentioned in the introduction in \eqref{Fs}.}
\begin{align}
 Ric^g &= f \cdot V^{\flat} \otimes V^{\flat}, \label{3} \\
  D \phi &= 0, \label{3bcd} \\
  V(f) &= 0, \label{4}
\end{align}
where $\text{Ric}=\text{Ric}^g$ denotes the Ricci tensor of $(M,g)$.
\end{Lemma}

\bprf
\eqref{3bcd} is a trivial consequence of $\phi$ being parallel and only listed for later reference. Every parallel spinor satisfies $\text{Ric}(X) \cdot \phi = 0$ (cf. \cite{ba81}). As also $V \cdot \phi = 0$, it follows that $g(V,Ric(X))=0$. Moreover, both $\text{Ric}(X)$ and $V$ are lightlike and orthogonal, and therefore they have to be linearly dependent. That is, there is $\theta \in \Omega^1(M)$ such that $\text{Ric}(X) = \theta(X) \cdot V$. In other words
\begin{align}
\text{Ric}(X,Y) = \theta(X) g(V,Y). \label{str}
\end{align}
Let $T$ be a timelike vector field on $M$ with $g(V,T) = 1$. \eqref{str} applied twice yields that $\theta(X) = \text{Ric}(X,T) = \text{Ric}(T,X) = \theta(T) g(V,X)$ which proves \eqref{3} for $f=\theta(T)$. Since $M$ admits a parallel spinor, we have that $\text{scal}^g = 0$. Thus\footnote{Here and in the following, $\delta = \delta^g$ denotes the divergence operator, i.e. given a $(k,0)$ tensor field $B$ on $(M,g)$, we set $(\delta B)(X_2,...,X_k) = - \sum_{a=0}^n \epsilon_a \left(\nabla^g_{s_a}B\right)(s_a,X_2,...,X_k)$ and for a vector field $V$ we have $\text{div}V = - \delta V^{\flat}$.}
\[ 0 = - (\delta (\text{Ric}))^{\sharp} \stackrel{\eqref{3}}{=} V(f)V + \text{div}V \cdot V + {\nabla}_V V. \] 
As $V$ is parallel, \eqref{4} follows.
\eprf
Let us now additionally assume that  $\Sigma \subset M$ is a spacelike hypersurface with induced Riemannian metric $g^{\Sigma}$ and future-directed unit normal vector field $T$ along $\Sigma$.  
For $X,Y\in T\Sigma$ denote by
\[
W(X,Y):=-g( \nabla_XT,Y)
\]
the second fundamental form of $(\Sigma,g^{\Sigma})\subset (M,g)$, i.e.,
we have
\[
\nabla_XY=\nabla^{\Sigma}_XY-W(X,Y) T,
\]
in which $\nabla^{\Sigma}$ denotes the Levi-Civita connection of $g^{\Sigma}$.
The dual of the  second fundamental form is the {\em Weingarten operator}, also denoted by $W$, and defined by
\[W(X,Y)=g^{\Sigma}(W(X),Y).
\]
It holds that $W=-\nabla T|_{T\Sigma}$.\\
\newline
Wrt. the spin structure canonically induced by that of $M$, let $(S^{\Sigma},\nabla^{\Sigma})$ be the spinor bundle of the
space-like hypersurface $(\Sigma,g^{\Sigma})$ with its spin derivative. Then
there is a canonical identification of $S^{\Sigma}$ with $ S_{|\Sigma}$
if $n$ is even and of $S^{\Sigma}$ with the half-spinors $
S^+_{|\Sigma}$ if $n$ is odd. For a detailed explanation of these standard identifications we refer
to \cite{bgm}. In this identification, the Clifford
product with a vector field $X$ on $\Sigma$ in both bundles is related
via
\begin{align*} 
 X \widetilde{\cdot} \varphi = i \, T \; {\cdot} \; X \; {\cdot} \;
\phi_{|\Sigma},
\end{align*}
where $\varphi \in \Gamma(S^{\Sigma})$ is
identified with $\phi_{|\Sigma} \in \Gamma( S^{(+)}_{|\Sigma})$. In
the following we will omit the $\widetilde{}$ subscript to the Clifford
multiplication in $\Sigma$ in order to keep the notation
simple, it will always be clear in which spinor bundle we are
working.  In general, the Dirac current $U_{\psi}$ of a spinor field $\psi$ on a
Riemannian spin manifold $(\M,\g)$ is given by
\begin{align*}
 \g(U_{\psi},X) := - i \,(X \cdot \psi,\psi) , \qquad X\in T\M.
 \end{align*}
In our situation, if $\phi \in \Gamma(S^{(+)})$ is a spinor field on $M$
and $\varphi:=\phi_{|\Sigma} \in \Gamma(S^{\Sigma})$ its restriction to $\Sigma$, the
Dirac currents are related by
\begin{align}
  (V_{\phi})_{|\Sigma} = \|\varphi\|^2 T_{|\Sigma} - U_{\varphi}. \label{vfg}
\end{align}

Using the above identification of the spinor bundles, the conditions
$\nabla \phi=0$ and $V_{\phi}\cdot \phi = 0$ translate into
the following conditions for the spinor field $\varphi= \phi_{|\Sigma}$ (cf. \cite{bam})

\begin{Proposition}
Let $(M,g)$ be a time-oriented Lorentzian spin manifold with parallel null spinor
field $\phi$. Then the spinor field $\varphi := \phi_{|\Sigma}$ on the
space-like hypersurface $(\Sigma,g^{\Sigma})$ satisfies
 \be \nabla^S_X
\varphi &=& \tfrac{i}{2}\, \W(X)\cdot \varphi \qquad
\forall \; X \in T\M, \label{spinor-1}\\
 U_{\varphi} \cdot \varphi &=& i\,u_{\varphi} \,\varphi, \label{spinor-2}
 \ee
where $\W$ is the Weingarten operator of $(\Sigma,g^{\Sigma})$ and $u_{\varphi}=
\sqrt{\g(U_{\varphi},U_{\varphi})} = \|\varphi\|^2$.
\end{Proposition}

We next evaluate \eqref{3} restricted to $\Sigma$ which gives $f_{|\Sigma}$: To this end, we use the following contracted versions of the Gau{\ss}- Codazzi- and Mainardi equation (cf. \cite{rin}).
\begin{Proposition}
Let $(\bM,\bg)$ be a time-oriented Lorentzian manifold with Einstein tensor $\overline{G} = \overline{Ric} - \frac{1}{2}\overline{scal} \cdot \bg$. Let $(\M,\g)$ be a spacelike hypersurface with induced metric, $T$ the future-directed timelike unit vector field along $\M$ and let $W$ be the second fundamental form of $(\M,g)$. Then on $\M$:
\begin{equation}
\begin{aligned} \label{hsf}
\overline{G}(T,T) & = \frac{1}{2} \left(\text{scal}^{\g} - \text{tr}_{\g}(W^2) +  (\text{tr}_{\g}W)^2 \right),\\
\overline{G}(T,X) & = (\delta^{\g} W) (X)  + d(\text{tr}_{\g}W)(X), \text{ }X \in T\M,
\end{aligned}
\end{equation}
where $\text{scal}^{\g}$ denotes the scalar curvature of $(\M,\g)$.
\end{Proposition}
We apply this result to the situation in Lemma \ref{eve} and combine it with equation \eqref{3}. As here $\text{scal}^{g}=0$ and therefore $G=\text{Ric}$, we find with $\text{scal}^{\Sigma}:=\text{scal}^{g^{\Sigma}}$ that 
\begin{align*}
  \frac{1}{2} \left(\text{scal}^{\Sigma} - \text{tr}_{g^{\Sigma}}(W^2) +  (\text{tr}_{g^{\Sigma}}W)^2 \right) = f_{|\Sigma} \cdot g(V_{|\Sigma},T)g(V_{|\Sigma},T) \stackrel{\eqref{vfg}}{=} f_{|\Sigma} \cdot u_{\ph}^2,
\end{align*}
that is
\begin{align}
f_{|\Sigma} =  \frac{1}{2u_{\ph}^2} \left(\text{scal}^{\Sigma} - \text{tr}_{g^{\Sigma}}(W^2) +  (\text{tr}_{g^{\Sigma}} W)^2 \right). \label{fsigma2}
\end{align}
Using the constraint equations only we can give another identity for $f_{|\Sigma}$ which becomes of importance later.
To this end we differentiate \eqref{spinor-1} again and skew-symmetrize to obtain
\[ R^{\Sigma}(X,Y) \ph = \frac{1}{4} \cdot (\W(X) \cdot \W(Y) - \W(Y) \cdot \W(X) ) \cdot \ph  + \frac{i}{2} \cdot \left((\nabla^{\Sigma}_X W)(Y) - (\nabla^{\Sigma}_Y W)(X) \right) \cdot \ph 
\]
for $X,Y \in T \Sigma$. Letting $Y=s_j$, taking the Clifford product with $s_j$ and summing over $j$ for some local orthonormal basis $(s_1,...,s_n)$ in $T\Sigma$ yields using that $W$ is $g^{\Sigma}$-symmetric that
\[ Ric^{\Sigma}(X) \cdot \ph = (W^2(X) \cdot \ph - \text{tr}_{g^{\Sigma}}(W) W(X) \cdot \ph) + i \cdot \left(X( \text{tr}_{g^{\Sigma}}(W)) \ph + \sum_{i=1}^n s_i \cdot (\nabla_{s_i}^{\Sigma} W)(X) \cdot \ph \right). \]
Taking another Clifford trace gives
\begin{align*}
-\text{scal}^{\Sigma} \cdot \ph = \left( (\text{tr}_{g^{\Sigma}}(W))^2 -  \text{tr}_{g_{\Sigma}}(W^2) \right) \cdot \ph + i \cdot \left(  2 d\text{tr}_{g^{\Sigma}}W  + 2 {\delta}^{\Sigma} W \right) \cdot \ph.
\end{align*}
Comparing with \eqref{fsigma2} gives
\begin{align*}
f_{|\Sigma} \cdot i \cdot {u_{\ph}^2} \cdot \ph = \left( d\text{tr}_{g^{\Sigma}}W  + {\delta}^{\Sigma} W \right) \cdot \ph
\end{align*}
Using the algebraic condition \eqref{spinor-2} on the left side, we arrive at
\begin{align*}
\left( d\text{tr}_{g^{\Sigma}}W  + {\delta}^{\Sigma} W - f_{|\Sigma} \cdot u_{\ph} \cdot U_{\ph}^{\flat} \right) \cdot \ph = 0.
\end{align*}
As $(\Sigma,g^{\Sigma})$ is Riemannian and $\ph$ nowhere vanishing, we conclude that 
\begin{align}
 d\text{tr}_{g^{\Sigma}}W  + {\delta}^{\Sigma} W  = f_{|\Sigma} \cdot u_{\ph} \cdot U^{\flat}_{\ph}. \label{leg}
\end{align}

\section{Proof of Theorem \ref{mother}}

The main ideas of the proof of Theorem \ref{mother} consist of the following steps:

\begin{enumerate}
 \item We impose the system of equations \eqref{3}-\eqref{4} as evolution equation for $(g,\phi,f)$ and rewrite it with the help of the background metric $h$ as a \textit{first} order quasilinear symmetric hyperbolic PDE of the form $A_0(t,x,u) \partial_t u = \sum_{\mu >0} A_{\mu}(t,x,u) \partial_{\mu} u + b(t,x,u)$, where $u$ collects the data $(g,\partial g, \phi,f)$. As initial data we choose $(g^{\Sigma},W,\ph,f_{|\Sigma})$ as defined or computed before. A local existence and uniqueness result for such PDEs is available in the smooth setting (cf. \cite{tay}) and gives locally defined $(g,\phi,f)$.
 \item We show that the locally given spinors $\phi$ are in fact parallel. This is achieved by showing that the data $\nabla_X \phi$ lie in the kernel of a locally defined second order normally hyperbolic (wrt. $g$) operator  $P$. Moreover, the generalized Killing spinor equation \eqref{spinor-1} shows that $\nabla \phi = 0$ initially. By a uniqueness result for the Cauchy problem for operators of type $P$ from \cite{bae}, we conclude that $\nabla \phi = 0$.
 \item Using the uniqueness of the local solutions obtained in (1), we can patch their domains together and obtain an open neighborhood $M$ of $\Sigma$ in $\mathbb{R} \times \Sigma$ on which $\phi$ is parallel. Furthermore, we can show that $M$ is globally hyperbolic. 
\end{enumerate}
Let us now carry out the technical details of each of these steps:\\
\newline
\textit{Step 1:}\\
The idea is to transform \eqref{3}-\eqref{4} into an evolution equation for the metric and spinor. However, there is a technical difficulty to overcome. The operator $g \mapsto \text{Ric}[g]=\text{Ric}^g$ mapping a metric to its Ricci-$(2,0)$-tensor is not easy to deal with from a PDE point of view: Locally, one has in coordinates that
\begin{align} \label{trustme}
\text{Ric}[g]_{\mu \nu} = - \frac{1}{2}g^{\alpha \beta}\partial_{\alpha} \partial_{\beta} g_{\mu \nu} + \nabla_{(\mu} \Gamma_{\nu)} + g^{\alpha \beta} g^{\gamma \delta} [ \Gamma_{\alpha \gamma \mu} \Gamma_{\beta \delta \nu} + \Gamma_{\alpha \gamma \mu} \Gamma_{\beta \nu \delta} + \Gamma_{\alpha \gamma \nu} \Gamma_{\beta \mu \delta} ],
\end{align}
where $\Gamma_{\nu}:=g^{\alpha \beta} \Gamma_{\alpha \nu \beta}$ and $\Gamma$ denote the Christoffel symbols of $g$ wrt. the coordinates.
From an analytic perspective, the fact that second order derivatives of $g$ appear in the summand $\nabla_{(\mu} \Gamma_{\nu)}$ prevents \eqref{3} considered as differential equation for $g$ in its present form from being accessible for hyperbolic PDE tools. To overcome this problem, we follow a strategy called \textit{hyperbolic reduction}, as applied in \cite{rin}, for instance, for the analogous problem for the Einstein field equations: To this end, we bring into play the fixed background metric
\begin{align}
  h:= -\lambda^2 dt^2 +  h_t. \label{hlam}
\end{align}
as in the formulation of the Theorem. Given local coordinates, we denote by $\widetilde{\Gamma}^{\mu}_{\alpha \beta}$ its Christoffel symbols. For any metric $g$ on $M$ we then introduce the difference tensor $A^{\mu}_{\alpha \beta} = \Gamma^{\mu}_{\alpha \beta} - \widetilde{\Gamma}^{\mu}_{\alpha \beta}$. Furthermore, we let 
\begin{align}
F_{\nu} &= g_{\mu \nu} g^{\alpha \beta} \widetilde{\Gamma}^{\mu}_{\alpha \beta}, \\
E_{\nu} &=-g_{\mu \nu}g^{\alpha \beta}A^{\mu}_{\alpha \beta}. \label{humm}
\end{align}
We denote by $\text{Sym}(\nabla E)[g]$ the symmetrization of the $(2,0)$-tensor $g(\nabla^g E, \cdot)$ for any given metric $g$. Then the operator
\[ \widehat{\text{Ric}}[g]:= \text{Ric}[g] + \text{Sym}(\nabla E)[g] \]
is by comparing with \eqref{trustme} in coordinates given by

\begin{align} 
 \widehat{\text{Ric}}_{\mu \nu} = - \frac{1}{2}g^{\alpha \beta}\partial_{\alpha} \partial_{\beta} g_{\mu \nu} + \nabla_{(\mu}F_{\nu)} + \underbrace{g^{\alpha \beta} g^{\gamma \delta} [ \Gamma_{\alpha \gamma \mu} \Gamma_{\beta \delta \nu} + \Gamma_{\alpha \gamma \mu} \Gamma_{\beta \nu \delta} + \Gamma_{\alpha \gamma \nu} \Gamma_{\beta \mu \delta} ]}_{=:H_{\mu \nu}[g,\partial g]}.
\end{align}
The crucial point is that second order derivatives of $g$ appear only in the first term of $\widehat{\text{Ric}}[g]$ (assured by addition of $E$; $F$ depends only on $g$ and not on its derivatives). Of course, eventually we will construct a solution and then show that $E=0$ as required by the Theorem.\\
\newline
With these technical preparations, we now show that under the assumptions of Theorem \ref{mother} every point $p \in \Sigma$ admits an open neighborhood $\mathcal{U}_p$ in $\mathbb{R} \times \Sigma$, a unique Lorentzian metric $g=g^{\mathcal{U}_p}$ , a spinor $\phi=\phi^{\mathcal{U}_p}$ and a function $f=f^{\mathcal{U}_p}$ defined on $\mathcal{U}_p$ subject to the PDE

\begin{align}
\text{Ric}^g  &= f \cdot V^{\flat} \otimes V^{\flat} - \text{Sym}(\nabla E),  \label{eq1}\\
{D}^g \phi & = 0, \label{eq2} \\
V(f) &=0, \label{eq3}
\end{align}

with initial data to be specified. To this end, choose $\mathcal{V}_p$ to be a coordinate neighborhood of $p \in \Sigma$ in $\mathbb{R} \times \Sigma$ with coordinates $(x_0=t,x_1,...,x_n)$, where $(x_1,...,x_n)$ are coordinates on $\Sigma$.
In the following, we have to distinguish various indices:
\begin{itemize}
\item Greek indices $\mu, \nu,...$ refer to the coordinates $(x_0,...,x_n)$, 
\item Latin indices $i,j,k,...$ refer to the spatial coordinates and run from 1 to $n$,
\item $a,b,c,...$ run from 0 to $n$ and appear as indices for local orthonormal bases.
\end{itemize}
 Given any metric $g$ sufficiently close\footnote{Concretely, the Gram Schmidt procedure wrt. $g$ should be applicable to the basis $(\partial_0,...,\partial_n)$.} to the fixed reference metric $h$, we form a $g$-dependent pseudo-orthonormal basis $(s_0,...,s_n)$ for $g$, i.e. $g(s_a,s_b) = \epsilon_a \delta_{ab}$, where $s_0$ is timelike and 
\begin{align}
  s_a = s_a[g] = \sum_{\mu} \zeta^{\mu}_a[g] \partial_{\mu} \label{rea}
\end{align}

on $\mathcal{V}_p$ where the defining coefficients $\zeta^{\mu}_i[g]$ can be chosen to depend smoothly on $g$ (via the Gram Schmidt procedure applied to $g$) and by the special form of the fixed background metric $h=-\lambda^2 dt^2 + h_t$ we may moreover assume that $\zeta^0_0[h]=\frac{1}{\lambda}$ and $\zeta^0_{a>0}[h]=0$. For such a metric, we rewrite \eqref{eq1}-\eqref{eq3} locally as follows:\\
\newline
To start with, \eqref{eq3} becomes equivalent to $V(f) = -\sum_a \epsilon_a \langle s_a \cdot \phi, \phi \rangle s_a(f) = 0$ which in terms of the coordinates can be rewritten as
\begin{align}
 - \sum_{a=0}^n \epsilon_a \zeta^0_a[g] \langle s_a \cdot \phi, \phi \rangle \partial_t f = \sum_{i}\sum_{a=0}^n \epsilon_a \zeta^{i}_a[g] \langle s_a \cdot \phi, \phi \rangle \partial_{i}f. \label{vf}
\end{align}
We next consider \eqref{eq2}, which in terms of the $s_a$ is equivalent to $\sum_{a=0}^n \epsilon_a s_a \cdot \nabla^g_{s_a} \phi = 0$. We trivialize the spinor bundle over $\mathcal{V}_p$ by means of the ONB $(s_0,...,s_n)$ and then think of $\phi$ wrt. this trivialization equivalently as a smooth function 
\[
\widetilde{\phi}=\widetilde{\phi}^{\mathcal{V}_p}: \mathcal{V}_p \rightarrow \Delta_{1,n} \cong \mathbb{C}^{2^{\frac{n}{2}}}.
\]
In this identification, we have that (cf. \cite{ba81}) $\nabla_{s_a} \phi = s_a(\widetilde{\phi}) + \omega(s_a) \cdot \widetilde{\phi}$, where the spin connection $\omega = \omega[g,\partial g]$ depends on the fixed metric and its first derivatives and $\omega(s_a) \in \mathbb{C}^{2^{\frac{n}{2}} \times 2^{\frac{n}{2}}}$, cf. \cite{ba81}. Let $\gamma_{a=0,...,n}$ be Gamma-matrices. That is, $\gamma_a \in \mathbb{C}^{2^{\frac{n}{2}} \times 2^{\frac{n}{2}}}$ and $\gamma_a \gamma_b + \gamma_b \gamma_a = -2 \epsilon_a \delta_{ab}$. Then $D \phi = 0$ is after multiplication with $s_0$ in the local trivialization equivalent to
\begin{align}
\sum_{a=0}^n \epsilon_a \zeta^0_a[g] \gamma_0 \cdot \gamma_a \cdot \partial_t \widetilde{\phi} = -\sum_{i} \sum_{a=0}^n \epsilon_a \zeta^{i}_a[g] \gamma_0 \cdot \gamma_a \cdot  \partial_{i} \widetilde{\phi} - \sum_{a=0}^n \epsilon_a \gamma_0 \cdot \gamma_a \cdot \omega(s_a) \cdot \widetilde{\phi}. \label{dph}
\end{align}
In order to rewrite \eqref{eq3} as system of first order equations, we introduce in analogy to \cite{first} for any Lorentzian metric $g$ on $\mathcal{V}_p$ close to $h$ the quantities $k_{\mu \nu}:= \partial_t g_{\mu \nu}$ and $g_{\mu \nu, i}:= \partial_{i} g_{\mu \nu}$, where $\rho$ runs only over the spatial indices $1,...,n$. In terms of these quantities, \eqref{eq1} can be rewritten as
\begin{align}
\partial_t g_{\mu \nu} &= k_{\mu \nu}, \label{f1} \\
g^{ij} \partial_t g_{\mu \nu, i} &= g^{ij} \partial_{i} k_{\mu \nu}, \label{f2} \\
- g^{00} \partial_t k_{\mu \nu} &= 2 g^{0 j} \partial_{j} k_{\mu \nu} + g^{ij} \partial_{j} g_{\mu \nu, i} - 2 H_{\mu \nu}[g,k] - 2 \nabla_{(\mu}F_{\nu)}[g,k] + 2f \cdot V_{\mu}V_{\nu}[g,\widetilde{\phi}], \label{f3}
\end{align}
and this system is equivalent to \eqref{eq1}. Indeed, let a triple $(g_{\mu \nu},k_{\mu \nu},g_{\mu\nu,i})$ solve \eqref{f1}-\eqref{f3}. As $g^{ij}$ is invertible for $g$ sufficiently close to $h$, \eqref{f2} is the same as $\partial_t g_{\mu \nu, i}=\partial_{i} k_{\mu \nu}$, and \eqref{f1} then gives $\partial_t (g_{\mu \nu, i}-\partial_{i} g_{\mu \nu})=0$. Appropriate choice of initial data (as done below) ensures $g_{\mu \nu, i}=\partial_{i} g_{\mu \nu}$ at $t=0$ and thus equality everywhere. Then \eqref{f3} is nothing but \eqref{eq1}. \\
\newline
The equations \eqref{vf}-\eqref{f3}, which are all defined on $\mathcal{V}_p$, can be summarized in the quasilinear first order PDE
\begin{align}
A_0(t,x,u) \partial_t u = \sum_{i = 1}^n A_{i}(t,x,u) \partial_{i} u + b(t,x,u), \label{system}
\end{align}
where we set 

\allowdisplaybreaks 
%\begin{equation}
\begin{align*} 
%\label{44}
u_{\mu \nu} &= \begin{pmatrix} g_{\mu \nu} \\ \begin{pmatrix} g_{\mu \nu, 1} \\ \vdots \\ g_{\mu \nu, n} \end{pmatrix} \\ k_{\mu \nu} \end{pmatrix}, \\  
u &= \begin{pmatrix} \begin{pmatrix} u_{00} \\ \vdots \\ u_{n+1 n+1} \end{pmatrix} \\ f \\ \widetilde{\phi} \end{pmatrix}, \\ 
\widetilde{A}_0(t,x,u) &=\begin{pmatrix} 1 & 0 & 0  \\ 0 & (g^{i j})_{i,j} & 0  \\ 0 & 0 & - g^{00}  \end{pmatrix}, \\
\widetilde{A}_i(t,x,u) &= \begin{pmatrix} 0 & 0 & 0 \\ 0 & 0 & \begin{pmatrix} g^{i1} \\ \vdots \\ g^{in} \end{pmatrix} \\ 0 & \begin{pmatrix} g^{i1} & \cdots & g^{in} \end{pmatrix} & 2 \cdot g^{0 i} \end{pmatrix}, \\
A_0(t,x,u) &= \begin{pmatrix} \widetilde{A}_0 \otimes Id_{(n+1)^2} & 0 & 0 \\ 0 & - \sum_{a=0}^n \epsilon_a \zeta^0_a[g] \langle s_a \cdot \phi, \phi \rangle & 0 \\ 0 & 0 & - \sum_{a=0}^n \epsilon_a \zeta^0_a[g]  \gamma_0 \cdot \gamma_a \end{pmatrix}, \\
A_{i}(t,x,u) &= \begin{pmatrix} \widetilde{A}_i \otimes Id_{(n+1)^2} & 0 & 0 \\ 0 &  \sum_{a=0}^n \zeta^{i}_a[g] \langle s_a \cdot \phi, \phi \rangle & 0 \\ 0 & 0 & \sum_{a=0}^n \epsilon_a \zeta^{i}_a[g] \gamma_0 \cdot \gamma_a \end{pmatrix},\\ 
b_{\mu \nu}(t,x,u) &= \begin{pmatrix} k_{\mu \nu} \\ 0 \\ - 2 H_{\mu \nu}[u] - 2 \nabla_{(\mu}F_{\nu)}[u] + 2f \cdot V_{\mu} V_{\nu}[u] \end{pmatrix}, \\
b(t,x,u) &= \begin{pmatrix} \begin{pmatrix} b_{00} \\ \vdots \\ b_{n+1 n+1} \end{pmatrix} \\  0  \\ \sum_{a=0}^n \epsilon_a \gamma_0 \cdot \gamma_a \cdot \omega(s_a)[u] \end{pmatrix}.
\end{align*}
%\end{equation}
%schematically 
%\begin{align}
%u &= \begin{pmatrix} g_{\mu \nu} \\ g_{\mu \nu, \rho} \\ k_{\mu \nu} \\ f \\ \widetilde{\ph} \end{pmatrix}, \\  
%A_0(t,x,u) &=\begin{pmatrix} I & 0 & 0 & 0 & 0 \\ 0 & g^{\rho \tau} & 0 & 0 & 0 \\ 0 & 0 & - \text{diag}(g^{00},...,g^{00}) & 0 & 0 \\ 0 & 0 & 0 & - \sum_{i=0}^n \epsilon_i a^0_i[g] \langle s_i \cdot \phi, \phi \rangle & 0 \\ 0 & 0 & 0 & 0 & - \sum_{i=0}^n \epsilon_i a^0_i[g] \langle s_i \cdot \phi, \phi \rangle  \end{pmatrix}, \\ \label{44}
%A_{\rho}(t,x,u) &= \begin{pmatrix} 0 & 0 & 0 & 0 & 0 \\ 0 & 0 & g^{\rho \tau} & 0 & 0 \\ 0 & g^{\rho \tau} & 2 \cdot \text{diag}(g^{0 \rho},...,g^{0 \rho}) & 0 & 0 \\ 0 & 0 & 0 & - \sum_{i=0}^n a^{\rho}_i[g] \langle s_i \cdot \phi, \phi \rangle & 0 \\ 0 & 0 & 0 & 0 & \sum_{i=0}^n \epsilon_i a^{\rho}_i[g] \gamma_0 \cdot \gamma_i \end{pmatrix},\\ \label{33}
%b(t,x,u) &= \begin{pmatrix} k_{\mu \nu} \\ 0 \\ - 2 H_{\mu \nu}[u] - 2 \nabla_{(\mu}F_{\nu)}[u] + 2f \cdot V_{\mu} V_{\nu}[u] \\ 0  \\ - \sum_{i=0}^n \epsilon_i \gamma_0 \cdot \gamma_i \cdot \omega(s_i)[u] \end{pmatrix},
%\end{align}
%where $I$ denotes the identity matrix in the respective dimension.\\
%\newline
We next specify initial data $u_0 := u_{|\Sigma \cap \mathcal{V}_p}$ for this first order PDE and then show that at least locally around the initial data the matrices $A_0$ and $A_{i}$ give rise to a symmetric hyperbolic PDE for which a solution theory is available.\\
\newline
\textit{Initial data:}
Concerning $g$ we set $g_{|\Sigma \cap \mathcal{V}_p} := h_{|\Sigma \cap \mathcal{V}_p}$, or in coordinates we set for $\lambda_{\Sigma}:=\lambda_{| \Sigma}$ 
\begin{align*}
{g_{ij}}_{|t=0} &= g^{\Sigma}(\partial_{i},\partial_{j}), \\
{g_{0 i}}_{|t=0} &= 0, \\
{g_{0 0}}_{|t=0} &= -\lambda_{\Sigma}^2.
\end{align*}
Motivated by the preceding discussion, we then set for $g_{\mu\nu, i}$
\[ {g_{\mu \nu, i}}_{|t=0}={\partial_{i} (g_{\mu \nu}}_{|t=0}). \]
Concerning $k_{\mu \nu}$ we observe that $\frac{1}{\lambda_{\Sigma}}\partial_t$ is the unit normal vector field wrt. $g$ along $\Sigma$ and set 
\begin{align}
{k_{ij}}_{|t=0} = -2 \lambda_{\Sigma}  W(\partial_i,\partial_j), \label{umms}
\end{align}
required, of course, by the fact that $(\Sigma,g^{\Sigma})$ should eventually embed with Weingarten tensor being the given $W$. The initial data for $k_{i 0}$ and $k_{00}$ are uniquely determined by the natural requirement ${(E_{\mu})}_{|t=0} \stackrel{!}{=} 0$  for any solution $g$. It is by definition of $E$ straightforward to compute that (cf. \cite{rin}) this is the case iff
\begin{equation}
\begin{aligned} \label{drop}
 {k_{00}}_{|t=0} & = -2 \lambda_{\Sigma}^2 \cdot {F_0}_{|t=0} + 2 \lambda_{\Sigma}^3 \cdot \text{tr}_{g^{\Sigma}} W, \\
 {k_{0 i}}_{|t=0} & = \lambda_{\Sigma}^2 \cdot \left[-F_{i} + \frac{1}{2} g^{jk} (2 \partial_{j} g_{ki} - \partial_{i} g_{jk}) +\partial_i(\text{log}\lambda_{\Sigma})\right]_{|t=0}.
\end{aligned}
\end{equation}
Note that it makes sense here to write $F_{|t=0}$ as the $F$-dependence on $g$ is only algebraic and $g_{|t=0}$ has already been specified. Finally, initial data for $f$ are motivated by \eqref{fsigma2}, i.e. we set 
\begin{align}
 f_{|t=0} =  \frac{1}{2u_{\ph}^2} \left(\text{scal}^{{\Sigma}} - \text{tr}_{g^{\Sigma}}(W^2) +  (\text{tr}_{g^{\Sigma}}W)^2 \right) \label{fsigma}.
\end{align}

\textit{Unique solvability of \eqref{system}}: The gamma matrices $\gamma_a \in \mathbb{C}^{2^\frac{n}{2} \times 2^\frac{n}{2}}$ can always be chosen such that $\gamma_0^{\dagger} = \gamma_0$ and $\gamma_a^{\dagger} = - \gamma_a$ for $a>0$. In fact, this is just a reformulation of $\langle x \cdot v, w \rangle = \langle v, x \cdot w \rangle$ (cf. \cite{ba81}) for the invariant inner product $\langle \cdot, \cdot \rangle$, where $v,w \in \Delta_{1,n}$ and $x \in \mathbb{R}^{1,n}$. It follows together with anticommutativity of the $\gamma_a$ that $(\gamma_0 \gamma_a)^{\dagger} = \gamma_0 \gamma_a$, whence the $\gamma_0 \cdot \gamma_a$ are symmetric when considered as real matrices. It follows that $A_0$ and $A_{i}$ as given above are symmetric matrices. Furthermore, we have with the initial conditions and recalling \eqref{rea} that\footnote{It is precisely this conclusion and its consequence \eqref{posded} which due to the signature of the metric do not hold in the analogous Riemannian situation with $f=0$.}
\begin{align*}
 A_0(t=0,x,u_0) = \begin{pmatrix} \begin{pmatrix} 1 & 0 & 0  \\ 0 & ({\left(g^{\Sigma}\right)}^{i j})_{i,j} & 0  \\ 0 & 0 & \lambda_{\Sigma}^2  \end{pmatrix} \otimes Id_{(n+1)^2} & 0 & 0  \\   0 & \frac{u_{\ph}}{\lambda_{\Sigma}} & 0 \\ 0 & 0 &  \frac{1}{\lambda_{\Sigma}} \cdot Id \end{pmatrix}.
\end{align*}
That is, we may restrict the fixed neighborhood $\mathcal{V}_p$ of $p$ in $\mathbb{R} \times \Sigma$ such that there is some constant $c$ with
\begin{align}
A_0(t,x,u) \geq cI > 0 \label{posded}
\end{align}
for all $(t,x) \in \mathcal{V}_p$ and $u$ sufficiently close to $u_0$. These observations, however, mean that \eqref{system} is a quasilinear symmetric hyperbolic system in the sense of \cite{tay}, Chapter 16.1-2, defined (by using charts) on a sufficiently small open subset of $\mathbb{R}^N$ for appropriate $N$ containing the initial data. By the existence and uniqueness result given in the reference, there is for smooth initial data $u_0$ around $p \in \Sigma$, as specified here, a neighborhood $\mathcal{U}_p \subset \mathcal{V}_p$ of $p$ in $\mathbb{R} \times \Sigma$ such that \eqref{system} has a unique smooth solution $u$ on $\mathcal{U}_p$ with $u_{|t=0} = u_0$. \\
Given this solution $u = \begin{pmatrix} g_{\mu \nu} & g_{\mu \nu, i} & k_{\mu \nu} & f & \widetilde{\ph} \end{pmatrix}$, we define with the coordinates $x_{\mu}$ on $\mathcal{V}_{p}$ specified earlier and with the ONB $s_a$ given by \eqref{rea} the bilinear from $g = g_{\mu \nu} dx^{\mu} dx^{\nu}$ on  $\mathcal{U}_p$. It is symmetric since \eqref{system} restricts to a PDE for symmetric matrices in the first 3 entries of $u$. Furthermore, after restricting $\mathcal{U}_p$ we may assume that $g$ is of Lorentzian signature on $\mathcal{U}_p$ as this holds for $g_{|\Sigma} = -\lambda_{\Sigma}^2 dt^2 + g^{\Sigma}$. After further restricting $\mathcal{U}_p$ to a smaller neighborhood of $p$, we may additionally assume that on $\mathcal{U}_p$
\begin{align}
g_{00} < 0, \\
(g_{ij})_{i,j > 0} > 0. \label{3sa}
\end{align}
We equip $\mathcal{U}_p$ with the spin structure $\mathcal{Q}_{\mathcal{U}_p}$ naturally induced by that of $\Sigma$ and trivialized by the $s_a$. On $\mathcal{U}_p$ we set
\begin{align*}
\phi = \phi^{\mathcal{U}_p} = [\widetilde{(s_0,...,s_n)},\widetilde{\phi}] \in \Gamma(\mathcal{U}_p,S^{g^{\mathcal{U}_p}}),
\end{align*}
where $\widetilde{(s_0,...,s_n)}$ denotes the lift of $(s_0,...,s_n)$ to $\mathcal{Q}_{\mathcal{U}_p}$.\\
\newline
Summarising the first step, we have constructed for each $p \in \Sigma$ an open neighborhood $\mathcal{U}_p$ of $p$ in $\mathbb{R} \times \Sigma$ and unique data $(g,f,\phi)$ defined on $\mathcal{U}_p$ solving \eqref{eq1}-\eqref{eq3} and which restrict to $(g^{\Sigma},\ph,f_{\Sigma})$ on $\Sigma$. Finally, it is possible to restrict $\mathcal{U}_p$ further (denoted by the same symbol) such that (cf. \cite{bae}) the spacelike hypersurface $\Sigma_p := \Sigma \cap \mathcal{U}_p$ is a Cauchy hypersurface in $(\mathcal{U}_p,g)$. By construction of the initial data \eqref{umms} and as $k_{\mu \nu} = \partial_t g_{\mu \nu}$, $(\Sigma_p,g^{\Sigma})$ embeds into $\mathcal{U}_p$ with Weingarten tensor $W$.\\
\newline
\textit{Step 2:}\\
In the terminology of the Step 1, we next show that the data $(\nabla \phi, V_{\phi} \cdot \phi, E)$ vanish on $\mathcal{U}_p$. All calculations and operators are wrt. the metric $g = g^{\mathcal{U}_p}$ on $\mathcal{U}_p$ as just specified. The idea is to construct a second order normally hyperbolic operator with the above data in its kernel. 
To this end, let $T:=s_0 := \frac{1}{\sqrt{-g(\partial_t,\partial_t)}}\partial_t$ be the unit timelike vector field on $\mathcal{U}_p$ which on $\Sigma_p$ coincides with $\frac{1}{\lambda_{\Sigma}}\partial_t$, the unit normal vector field to $\Sigma_p$ wrt. $g$. We trivialize $T^{\bot}$ via an orthonormal basis $(s_1,...,s_n)$ and it follows that at $t=0$ the $(s_1,...,s_n)$ are a pointwise ONB for $(\Sigma_p,g^{\Sigma})$. \\
We start with finding a second order equation for $E$. Let $G = \text{Ric}- \frac{1}{2}\text{scal} \cdot g$ denote the Einstein tensor of $(\mathcal{U}_p,g)$. 
By the first line of \eqref{eq1} we get for $X, Y \in T\mathcal{U}_p$
\begin{equation}
\begin{aligned} \label{treu}
G(X,Y) - f g(V,X)g(V,Y) + \frac{1}{2}f g(V,V) \cdot g(X,Y) =& -\frac{1}{2} g(\nabla_{X}E,Y) -\frac{1}{2} g(\nabla_{Y}E,X)\\
&+ \frac{1}{2} \text{tr}_g(\nabla E) \cdot g(X,Y), 
\end{aligned}
\end{equation}
where we abbreviate $V:=V_{\phi}$. Let us from now on denote $E$ and its metric dual with the same symbol. We take the divergence of the $(2,0)$-tensors in \eqref{treu} on both sides, use that $\delta G = 0$, to obtain
\begin{align}
-V(f) V  - f \cdot \text{div} V \cdot V - f \cdot \nabla_V V + \frac{1}{2} \text{grad}f \cdot g(V,V) + f \cdot \text{grad}||V||^2 =  \frac{1}{2} \Delta^{\nabla} E - \frac{1}{2} \text{Ric}(E), \label{E1}
\end{align}
where $\Delta^{\nabla} = - \sum_a \epsilon_a \left(\nabla_{s_a} \nabla_{s_a} + \text{div}(s_a) \cdot \nabla_{s_a} \right)$ denotes the Bochner Laplacian. However, as a simple consequence of $D \phi = 0$, we get that $\text{div} V_{\phi} = -2 \text{Re} \langle D\phi, \phi \rangle  = 0$. Moreover,
\begin{align*}
\nabla_{s_d}V_{\phi} = -2 \cdot \sum_{a=0}^n \epsilon_a \text{Re} \langle s_a \cdot \nabla_{s_d} \phi, \phi \rangle \cdot s_a.
\end{align*}
$V(f)=0$ by construction in step 1. This allows rewriting \eqref{E1} as 
%Moreover, we use the Weizenboeck formula $ \Delta_1 = (d + \delta)^2 = \Delta^{\nabla} + Ric$ and $V(f) = 0$ (cf...) to conclude that 
%\begin{align}
\allowdisplaybreaks

\begin{equation}
\begin{aligned} \label{enab}
\Delta^{\nabla} E &=  \text{Ric}(E) - 2f \cdot \nabla_V V + \text{grad} f \cdot g(V,V) + 2f \cdot \text{grad}||V||^2 \\
 &= \text{Ric}(E) +4 f \cdot \sum_{a=0}^n \epsilon_a \text{Re} \langle s_a \cdot \nabla_{V} \phi, \phi \rangle \cdot s_a + \text{grad} f \cdot g(V,V) - 8f \cdot \sum_{a=0}^n \epsilon_a \text{Re} \langle V \cdot \nabla_{s_a} \phi, \phi \rangle s_a. 
\end{aligned}
\end{equation}
Using this, it is straightforward to compute for $d=0,...,n$:

\begin{equation}
\begin{aligned}
\Delta^{\nabla} (\nabla_{s_d} E)  =& \nabla_{s_d} (\Delta^{\nabla} E) + \sum_{a=0}^n \epsilon_a (R(s_d,s_a) \nabla_{s_a}E  - \nabla_{[s_a,s_d]} \nabla_{s_a} E - \nabla_{s_a} \nabla_{[s_a,s_d]}E  \\
&+ \nabla_{s_a} (R(s_d,s_a)E) + s_d(\text{div} s_a) \cdot \nabla_{s_a}E - \text{div}(s_a) \cdot \left(R(s_a,s_d)E + \nabla_{[s_a,s_d]}E \right)) \\
\stackrel{\eqref{enab}}{=}& \nabla_{s_d} \left(\text{Ric}(E) +4 f \cdot \sum_{a=0}^n \epsilon_i \text{Re} \langle s_a \cdot \nabla_{V} \phi, \phi \rangle \cdot s_a + \text{grad} f \cdot g(V,V) - 8f \cdot \sum_{a=0}^n \epsilon_a \text{Re} \langle V \cdot \nabla_{s_a} \phi, \phi \rangle s_a  \right) \\
& + \sum_{a=0}^n \epsilon_a \left( R(s_d,s_a) \nabla_{s_a}E  - \nabla_{[s_a,s_d]} \nabla_{s_a} E - \nabla_{s_a} \nabla_{[s_a,s_d]}E  + \nabla_{s_a} (R(s_d,s_a)E) + s_d(\text{div} s_a) \cdot \nabla_{s_a}E \right. \\
& \left. - \text{div}(s_a) \cdot \left(R(s_a,s_d)E + \nabla_{[s_a,s_d]}E \right) \right)
\end{aligned}
\end{equation}
We now turn to $\nabla \phi$ and compute for $d \in \{0,...,n\}$:
\begin{align*}
{D} \left(\nabla_{s_d} \phi \right) & = \sum_{a=0}^n \epsilon_a s_a \cdot \nabla_{s_a} \nabla_{s_d} \phi \\
& = \sum_{a=0}^n \epsilon_a s_a \cdot \nabla_{s_d} \nabla_{s_a} \phi + \sum_{a=0}^n \epsilon_a s_a \cdot R^S(s_a,s_d) \phi + \sum_{a=0}^n \epsilon_a s_a \cdot \nabla_{[s_a,s_d]} \phi  \\
& = \nabla_{s_d} \left({D} \phi \right) - \sum_{a=0}^n \epsilon_a \left(\nabla_{s_d} s_a \right) \cdot  \nabla_{s_a} \phi + \frac{1}{2} \cdot \text{Ric}(s_d) \cdot \phi + \sum_{a,b=0}^n \epsilon_a \epsilon_b g([s_a,s_d],s_b) s_a \cdot \nabla_{s_b} \phi. \label{eq1}
\end{align*}
Using \eqref{eq1}, we get that
\[ \text{Ric}(s_d) \cdot \phi = f \cdot g(V,s_d) \cdot V \cdot \phi - \frac{1}{2} (\nabla_{s_d}E)\cdot \phi - \frac{1}{2} \sum_{a=0}^n \epsilon_a g(\nabla_{s_a}E,s_d) s_a \cdot \phi. \] 
Moreover, we use $D \phi = 0$. This yields
\begin {equation}
\begin{aligned}
D^2 \left(\nabla_{s_d} \phi \right) =& - \sum_{a=0}^n \epsilon_a D\left( \left(\nabla_{s_d} s_a \right) \cdot  \nabla_{s_a} \phi \right) + \frac{1}{2} \cdot D \left( f \cdot g(V,s_d) \cdot V \cdot \phi - \frac{1}{2} (\nabla_{s_d}E)\cdot \phi - \frac{1}{2} \sum_{a=0}^n \epsilon_a g(\nabla_{s_a}E,s_d) s_a \cdot \phi \right) \\
& + \sum_{a,b=0}^n \epsilon_a \epsilon_b D \left( g([s_a,s_d],s_b) s_a \cdot \nabla_{s_b} \phi \right).
\end{aligned}
\end{equation}

As a next step we consider the spinor $V \cdot \phi$. %V = V phi
We compute
\begin{align*}
{D} (V \cdot \phi) &= \sum_{a=0}^n \epsilon_a \left(s_a \cdot \left(\nabla_{s_a} V \right) \cdot \phi + s_a \cdot V \cdot \nabla_{s_a} \phi\right) \\
& \stackrel{D \phi = 0}{=} -2 \sum_{a,b=0}^n \epsilon_a \epsilon_b \text{Re} \langle s_b \cdot \nabla_{s_a} \phi, \phi \rangle s_a \cdot s_b \cdot \phi - 2 \cdot \nabla_{V} \phi.
%&= \sum_{a=0}^n \epsilon_i \left( \left(s_i \wedge \nabla^{g}_{s_i} V^{\flat} \right) - s_i \invneg \nabla^{g}_{s_i} V^{\flat} \right) \cdot \phi -V \cdot s_i \cdot \nabla_{s_i} \phi -2 g(s_i,V) \cdot \nabla_{s_i} \phi \\
%&= \left( (d+ \delta)V^{\flat} \right) \cdot \phi - V \cdot {D} \phi - 2 \nabla_V \phi. 
%& = 
\end{align*}
Applying $D$ again thus gives 
\begin{align}
D^2(V \cdot \phi) &= -2 \sum_{i,j=0}^n \epsilon_a \epsilon_b D\left((\text{Re} \langle s_b \cdot \nabla_{s_a} \phi, \phi \rangle s_a \cdot s_b \cdot \phi \right) -2  \sum_{a=0}^n \epsilon_a D (g(s_a,V) \cdot \nabla_{s_a} \phi).
\end{align}
%Using this and $D \ph = 0$, we compute
%\begin{align}
%D(\nabla_{s_d}E \cdot \ph) & = \left( (d+ \delta)\nabla_{s_d}E \right) \cdot \phi - 2 \nabla_{\nabla_{s_d}E}^{{S}} \phi \\
%& = \sum_{i} \left(s_a^{\flat} \wedge \nabla_{s_a} \nabla_{s_d} E + s_a \invneg \nabla_{s_a} \nabla_{s_d} E \right) \cdot \phi - 2 \nabla_{\nabla_{s_d}E}^{{S}} \phi \\
%& = \left(\nabla_{s_b} dE + \sum_a s_a \wedge \left(R(X,s_a)E + \nabla_{\nabla_{s_a} s_b} E \right) + \nabla_{s_b} \delta E + Ric(s_b,E) + g(\nabla_{\nabla_{s_a} s_b} E,s_a) \right) \cdot \ph \\
%&- 2 \nabla_{\nabla_{s_d}E}^{{S}} \phi
%\end{align}
%In similar fashion one obtains
%\begin{align}
%D(g(\nabla_{s_a}E,s_d) s_a \cdot \phi) &= \left( (d+ \delta)g(\nabla_{s_a}E,s_d) s_a \right) \cdot \phi - 2 \nabla_{g(\nabla_{s_a}E,s_d) s_a}^{{S}} \phi \\
 %&= \left(R(E,X) + \sum_{i,j} g(\nabla_{s_a}E,\nabla_{s_b}s_d) s_a \wedge s_b + g(\Delta^{\nabla}E,s_b) + \sum_a g(\nabla_{s_a}E,\nabla_{s_a} s_b) \right) \cdot \ph
%\end{align}

Finally, we consider the function $||V||^2 = g(V,V)$. We have for the Laplacian $\Delta$ acting on functions
\begin{equation}
\begin{aligned} \label{deff}
\Delta ||V||^2 & = - \text{div}(\text{grad}(g(V,V))) \\
&= 4 \cdot \sum_{a=0}^n \epsilon_a \cdot \text{div} \left( \text{Re} \langle V \cdot \nabla_{s_a} \phi, \phi \rangle \cdot s_a \right). 
\end{aligned}
\end{equation}

In terms of differential operators, we summarize our previous computations as follows: Let us introduce for $d=0,...,n$ the notation 
\begin{align*}
\alpha_d &:= \nabla_{s_d} \phi, \\
\beta_d &:= \nabla_{s_d} E, \\
\chi &:= V_{\phi} \cdot \phi, \\
\eta &:= E, \\
\kappa &:= g(V_{\phi},V_{\phi}).
\end{align*}
In terms of these data the equations \eqref{enab}-\eqref{deff} can be rewritten as

\begin{equation}
\begin{aligned} \label{mastere}
0 =& D^2 \alpha_d + \sum_{a=0}^n \epsilon_a D\left( \left(\nabla^{{S}}_{s_d} s_a \right) \cdot  \alpha_a \right) - \frac{1}{2} \cdot D \left( f \cdot g(V,s_d) \cdot \chi -\frac{1}{2} \cdot \beta_d \cdot \phi - \frac{1}{2} \cdot \sum_{a=0}^n \epsilon_a g(\beta_a ,s_d) s_a \cdot \phi \right) \\
& - \sum_{a,b=0}^n \epsilon_a \epsilon_b D \left( g([s_a,s_d],s_b) s_a \cdot \alpha_b \right), \\
0 =& \Delta^{\nabla}\beta_d - \nabla_{s_d} \left(\text{Ric}(\eta) +4 f \cdot \sum_{a,b=0}^n \epsilon_a \epsilon_b g(V,s_b) \text{Re} \langle s_a \cdot \alpha_b, \phi \rangle \cdot s_a + \text{grad} f \cdot \kappa - 8f \cdot \sum_{a=0}^n \epsilon_a \text{Re} \langle V \cdot \alpha_a, \phi \rangle s_a  \right) \\
& - \sum_{a=0}^n \epsilon_a \left( R(s_d,s_a) \beta_a  - \nabla_{[s_a,s_d]} \beta_a - \nabla_{s_a}\left(\sum_{c=0}^n \epsilon_c g([s_a,s_d],s_c) \beta_c \right) + \nabla_{s_a} (R(s_d,s_a)\eta) + s_d(\text{div}(s_a)) \cdot \beta_a \right. \\
& \left. - \text{div}(s_a) \cdot \left(R(s_a,s_d)\eta + \nabla_{[s_a,s_d]} \eta \right) \right), \\
0 = & D^2 \chi + 2 \sum_{a,b=0}^n \epsilon_a \epsilon_b D\left(\text{Re} \langle s_b \cdot \alpha_a, \phi \rangle s_a \cdot s_b \cdot \phi \right) +2  \sum_{a=0}^n \epsilon_a D (g(s_a,V) \cdot \alpha_a), \\
0 = & \Delta^{\nabla} \eta -  \text{Ric}(\eta) -4 f \cdot \sum_{a,b=0}^n \epsilon_a \epsilon_b g(V,s_b) \text{Re} \langle s_a \cdot \alpha_b, \phi \rangle \cdot s_a - \kappa \cdot \text{grad}f  + 8f \cdot \sum_{a=0}^n \epsilon_a \text{Re} \langle V \cdot \alpha_a, \phi \rangle s_a, \\
0 =& \Delta \kappa - 4 \sum_{a=0}^n \epsilon_a \text{div} \left( \text{Re} \langle V \cdot \alpha_a, \phi \rangle \cdot s_a \right). 
\end{aligned}
\end{equation}
Another way of looking at this system goes as follows: We introduce the vector bundle 
\[ {\mathcal{E}}_{\mathcal{U}_p} := \left(\oplus_{d=0}^n S_{\mathcal{U}_p} \right) \oplus \left(\oplus_{d=0}^n T\mathcal{U}_p \right) \oplus S_{\mathcal{U}_p} \oplus T\mathcal{U}_p \oplus \underline{\mathbb{R}} \rightarrow \mathcal{U}_p. \]
It carries a covariant derivative naturally induced by $\nabla^g$ (and $\nabla^S$) and denoted by the same symbol. Let now $v=((\alpha_d)_{d=0,...,n},(\beta_d)_{d=0,...,n},\chi, \eta, \kappa)$ be an \textit{arbitrary} section of ${\mathcal{E}}_{\mathcal{U}_p}$. Then the right side of \eqref{mastere} can be interpreted as $P_{\mathcal{U}_p}v$, meaning the action of a second order \textit{linear} differential operator, 
\[ P_{\mathcal{U}_p} : \Gamma(\mathcal{U}_p,{\mathcal{E}}_{\mathcal{U}_p}) \rightarrow \Gamma(\mathcal{U}_p,{\mathcal{E}}_{\mathcal{U}_p}), \]
 on $v$, which wrt. the decomposition of ${\mathcal{E}}_{\mathcal{U}_p}$ into summands is given by
 \begin{align} \label{pdef}
 P_{\mathcal{U}_p} = \begin{pmatrix} D^2 & & & & & & & &\\ & \ddots & & & & & & & \\ & & D^2 & & & & & & \\& & & \Delta^{\nabla} & & & & & \\ & & & & \ddots & & & & \\ & & & & & \Delta^{\nabla} & &  & \\ & & & & & & D^2  & & \\ & & & & & & & \Delta^{\nabla} &  \\& & & & & & & & \Delta \end{pmatrix} + P^1_{\mathcal{U}_p},
 \end{align}
where $P^1_{\mathcal{U}_p} : \Gamma(\mathcal{U}_p,{\mathcal{E}}_{\mathcal{U}_p}) \rightarrow \Gamma(\mathcal{U}_p,{\mathcal{E}}_{\mathcal{U}_p})$ is a linear differential operator of order 1 whose concrete form can be extracted from \eqref{mastere} but it is irrelevant in the following. Then a reformulation of \eqref{enab}-\eqref{deff} is that the section 
\begin{align} \label{secu}
u= \begin{pmatrix} (\nabla_{s_d} \phi)_{d=0,...,n} \\ (\nabla_{s_d}E)_{d=0,...,n} \\ V \cdot \phi \\  E \\ g(V,V) \end{pmatrix} \in \Gamma(\mathcal{U}_p, {\mathcal{E}}_{\mathcal{U}_p})
\end{align}
 satisfies
\begin{align}
P_{\mathcal{U}_p} u = 0.
\end{align}
A crucial property of $P_{\mathcal{U}_p}$ defined above is that it is \textit{normally hyperbolic}. In general, given any Lorentzian manifold $(M,g)$ with smooth vector bundle ${\mathcal{E}} \rightarrow M$ and second order differential operator $P: \Gamma(\mathcal{E}) \rightarrow \Gamma(\mathcal{E})$, we say that $P$ is normally hyperbolic if its principal symbol is given by the metric, $\sigma_P(\zeta_x) = -g(\zeta_x,\zeta_x) \cdot id_{\mathcal{E}_x}$ for $\zeta_x \in T^*_xM$, or equivalently, if in local coordinates $(x^0,...,x^n)$ on $M$ and a local trivialization of $\mathcal{E}$ we have
\[ P = - \sum_{\mu, \nu = 0}^{n} g^{\mu \nu}(x) \frac{\partial^2}{\partial_{\mu} \partial_{\nu}} + \sum_{\mu=0}^{n} A_{\mu}(x) \frac{\partial}{\partial x^{\mu}} + B(x) \]
for matrix-valued coefficients $A_j$ and $B$ depending smoothly on $x$. It is well known that $D^2$ acting on spinors as well as $\Delta^{\nabla}$ acting on vector fields and $\Delta$ acting on functions enjoy this property. Normal hyperbolicity is moreover preserved when adding a differential operator of order at most 1. From \eqref{pdef} it then follows immediately that $P_{\mathcal{U}_p}$ is normally hyperbolic on $(\mathcal{U}_p,g)$ as well.\\
\newline
\textit{Initial data.} In order to apply a uniqueness result for solutions to $P_{\mathcal{U}_p} u = 0$ we next show that the section \eqref{secu} satisfies 
\begin{align}
u_{|\Sigma_p} &= 0, \label{f12} \\
\left(\nabla_T u \right)_{|\Sigma_p} &=0.  \label{f13}
\end{align}
As a consequence of the initial data imposed on $g$ in the first step, $\Sigma_p = \Sigma \cap \mathcal{U}_p$ embeds into $\mathcal{U}_p$ with Weingarten tensor $W$. But then $(\nabla_{s_a} \phi)_{|\Sigma_p} = 0$ follows for $a>0$ as this is just a reformulation of the generalized imaginary Killing spinor condition \eqref{spinor-1} for $\ph$ in terms of data on $(\mathcal{U}_p,g)$. Moreover, as $D \phi = 0$ on $\mathcal{U}_p$ by the first step, we have that $0 = (D \phi)_{|\Sigma_p} = - (T \cdot \nabla_T \phi)_{|\Sigma_p}$, i.e. also $(\nabla_{s_0} \phi)_{|\Sigma_p} = 0$. The algebraic constraint \eqref{spinor-2} rewritten in terms of $(M,g)$ precisely yields that $V_{\phi} \cdot \phi = 0$ on $\Sigma_p$. Multiplying this by $V_{\phi}$ again gives that also $g(V_{\phi},V_{\phi})(\Sigma_p) = 0$. We turn to the $E-$terms. The initial conditions for $\partial_t g_{00}$ and $\partial_t g_{0 j}$ (cf. \eqref{drop}) where chosen in such a way that $E_{|\Sigma_p}$ follows (insert \eqref{drop} into the definition of $E$ in \eqref{humm}). From this it already follows that $\nabla_{s_a}E$ vanishes on $\Sigma_p$ for $a>0$. The most involved part is to show that $\nabla_{s_0}E$ vanishes on $\Sigma_p$ as well: The data $(g,\phi,f)$ constructed in the first step solve \eqref{treu}. Inserting $T_{|\Sigma_p}=\frac{1}{\lambda_{\Sigma}}(\partial_t)_{|\Sigma_p}$ and $X \in T\Sigma_p$ into this $(2,0)$-tensor and using the hypersurface formulas \eqref{hsf} as well as $u_{\ph} = - g(V,T)_{|\Sigma}$ gives 

\begin{align} \label{exp}
(d\text{tr}_{g^{\Sigma}}W) (X)  + ({\delta}^{\Sigma} W)(X) - f_{|\Sigma} \cdot u_{\ph} \cdot g^{\Sigma}(U_{\ph},X) = -\frac{1}{2} h(\nabla_{X}E,T) -\frac{1}{2} h(\nabla_{T}E,X).
\end{align}
The left hand side of \eqref{exp} is zero due to \eqref{leg} and for the right hand side we obtain using $E_{|\Sigma_p} = 0$ as well as $(\nabla_XE)_{|\Sigma_p} = 0$
that 
\begin{align}
h(\nabla_{T}E,X) = 0. \label{67}
\end{align}
Similarly, inserting $(T_{|\Sigma_p}, T_{|\Sigma_p})$ into \eqref{treu} and evaluating on $\Sigma_p$ yields
\begin{align*}
\frac{1}{2} \left(\text{scal}^{g^{\Sigma}} - \text{tr}_{g^{\Sigma}}(W^2) +  (\text{tr}_{g^{\Sigma}}W)^2 \right) - f_{|\Sigma} \cdot u_{\ph}^2  = - h(\nabla_{T}E,T) - \frac{1}{2} \text{tr}_h(\nabla E). 
\end{align*}
The left side vanishes as this was precisely the initial condition for $f$ (cf. \eqref{fsigma2}). As $\nabla_{s_a}E = 0$ for $a>0$ the second summand on the right side reduces to $+\frac{1}{2} h(\nabla_{T}E,T)$, and in total
\[ h(\nabla_{T}E,T) = 0 \]
Combined with \eqref{67} this yields $\nabla_{T}E = 0$ on $\Sigma_p$. These observations prove \eqref{f12}.\\
\newline
We turn to $\nabla_T u$ on $\Sigma_p$. Concerning $\nabla_{s_d} \phi$, we find for $d > 0$ that 
\begin{align*}
\left( \nabla_T \nabla_{s_d} \phi \right)_{|\Sigma_p} & = \left( \nabla_{s_d} \nabla_{T} \phi + R^{S}(T,s_d) \phi + \nabla_{[T,s_d]} \phi \right)_{|\Sigma_p} \\
& = R^S(T,s_d) \phi_{|\Sigma_p},
\end{align*}
where the last equality follows from $\nabla \phi = 0$ on $\Sigma_p$. The remaining curvature term is evaluated as follows: \eqref{eq1} evaluated on $\Sigma_p$ yields with $u_{| \Sigma_p} = 0$ that 
\[ 0 = \text{Ric}(s_d) \cdot \phi = - 2 \cdot \sum_{a=0}^n \epsilon_a s_a \cdot R^S(s_d,s_a) \phi\]
on $\Sigma_p$.  It follows from $T^2 = 1$ that on $\Sigma_p$
\begin{align*}
 R^S(T,s_d) \phi = T \cdot \sum_{a=1}^n s_a \cdot R^S(s_a,s_d) \phi.
\end{align*}
However, every summand on the right hand side vanishes on $\Sigma_p$ as follows from differentiating $\nabla_{s_a} \phi$, which vanishes on $\Sigma_p$, in direction of $\Sigma_p$ again and skew-symmetrizing. Thus, also $R(T,s_d) \phi = 0$ on $\Sigma_p$ and therefore 
\begin{align}
\left( \nabla_T \nabla_{s_{d>0}} \phi \right)_{|\Sigma_p} = 0. \label{truef}
\end{align}

Moreover, tracing the equation $(\text{Ric}(X) \cdot \phi)_{|\Sigma_p} = 0$, which holds for every $X \in \mathfrak{X}({\mathcal{U}_p})$ yields that $\text{scal} \cdot \phi_{|\Sigma_p} = 0$ and it then follows with the Schr{\"o}der Lichnerowicz formula that on $\Sigma_p$
\begin{align}
0 = {D}^2 \phi = \Delta \phi = -\sum_{a=0}^n \epsilon_a \left( \nabla_{s_a} \nabla_{s_a} \phi + \text{div}(s_a) \cdot \nabla_{s_a} \phi \right). \label{a2}
\end{align}
On $\Sigma_p$, however, we have that $\nabla_{s_a} \nabla_{s_a} \phi = 0$ for $a > 0$ and $\nabla_{s_a} \phi = 0$ for all $a$. Thus, \eqref{a2} yields $\nabla_{T} \nabla_{T} \phi = 0$ on $\Sigma_p$. \\
Let us turn to the $E$-terms. It has already been shown that $(\nabla E)_{|\Sigma_p} = 0$. We get for $d>0$ that evaluated at $\Sigma_p$
\begin{align*}
\nabla_T \nabla_{s_d} E &= \nabla_{s_d} \nabla_T E + R(T,s_d)E + \nabla_{[T,s_d]}E \\
& = 0,
\end{align*}
as also $E_{|\Sigma_p } = 0$. Moreover, we have by simply rewriting \eqref{enab} that 
\begin{align*}
\nabla_T \nabla_T E =& \sum_{a=1}^n \epsilon_a \left(\nabla_{s_a} \nabla_{s_a} E + \text{div}(s_a) \cdot \nabla_{s_a}E \right) - \text{div}(T) \cdot \nabla_{T} E + \text{Ric}(E) +4 f \cdot \sum_{a=0}^n \epsilon_a \text{Re} \langle s_a \cdot \nabla_{V} \phi, \phi \rangle \cdot s_a \\
 &+ \text{grad} f \cdot g(V,V) - 8f \cdot \sum_{a=0}^n \epsilon_a \text{Re} \langle V \cdot \nabla_{s_a} \phi, \phi \rangle s_a. 
\end{align*}
However, all the terms on the right hand side vanish on $\Sigma_p$ due to the initial conditions that were already verified. Finally, $\nabla_T(V \cdot \phi)$ and $T(g(V,V))$ vanish on $\Sigma_p$ as these quantities are just linear expressions in $\nabla_T \phi$. This shows \eqref{f13}.\\
\newline
\textit{Cauchy problem for }$P_{\mathcal{U}_p}:$ We have shown that the operator $P_{\mathcal{U}_p}$ and the section $u$ defined by \eqref{secu} satisfy
\begin{equation}
\begin{aligned} \label{34}
P_{\mathcal{U}_p}u &= 0, \\
u_{| \Sigma_p} & = 0, \\
\left(\nabla_T u \right)_{|\Sigma_p} &=0.
\end{aligned}
\end{equation}
Moreover, $(\mathcal{U}_p,g)$ was constructed to be globally hyperbolic with spacelike Cauchy hypersurface $\Sigma_p$. The main result of \cite{bae} states that for this geometry the Cauchy problem for the normally hyperbolic operator $P_{\mathcal{U}_p}$ is well-posed, i.e. the system \eqref{34} admits a unique solution. As $P_{\mathcal{U}_p}$ is linear, it is clear that the zero section solves \eqref{34}. As the solution is unique, we conclude that 
\begin{align} 
u= \begin{pmatrix} (\nabla_{s_d} \phi)_{d=0,...,n} \\ (\nabla_{s_d}E)_{d=0,...,n} \\ V \cdot \phi \\  E \\ g(V,V) \end{pmatrix} \equiv 0
\end{align}
on $\mathcal{U}_p$. In particular, $\nabla \phi = 0$ on $\mathcal{U}_p$ with lightlike Dirac current.\\
\newline
\textit{Step 3:}\\
We next globalize the local development of the initial data. In step 1 we have constructed for every $p \in \Sigma$ the data $(g^{\mathcal{U}_p},\phi^{\mathcal{U}_p},f^{\mathcal{U}_p})$ defined on some open $\mathcal{U}_p \subset \mathbb{R} \times \Sigma$ sufficiently small. Let $p,q \in \Sigma$ and assume that $\mathcal{U}_p \cap \mathcal{U}_q \neq \emptyset$. Choose coordinates $(x_0,...,x_n)$ and $(y_0,...,y_n)$ on $\mathcal{U}_p$ and $\mathcal{U}_q$ respectively as in step 1. On $\mathcal{U}_p \cap \mathcal{U}_q$ we consider the coordinates given by restriction of the $x_i$. Then wrt. these coordinates the data
\begin{align}
u_p = \begin{pmatrix} g^{\mathcal{U}_p}_{\mu \nu} \\ g^{\mathcal{U}_p}_{\mu \nu, i} \\ k^{\mathcal{U}_p}_{\mu \nu} \\ f^{\mathcal{U}_p} \\ \widetilde{\phi}^{\mathcal{U}_p} \end{pmatrix}, \text{ }
u_q = \begin{pmatrix} g^{\mathcal{U}_q}_{\mu \nu} \\ g^{\mathcal{U}_q}_{\mu \nu, i} \\ k^{\mathcal{U}_q}_{\mu \nu} \\ f^{\mathcal{U}_q} \\ \widetilde{\phi}^{\mathcal{U}_q} \end{pmatrix}
\end{align}
as introduced in step 1 solve by construction the system \eqref{system} formulated in the $x-$coordinates. This follows as the system \eqref{system} is manifestly coordinate invariant as it can also be equivalently formulated in terms of \eqref{eq1}-\eqref{eq3}. Moreover the initial data $(u_p)_{|\Sigma_p} = (u_q)_{|\Sigma_q}$ coincide since they arise as restrictions of globally defined data on $\Sigma$. It then follows from the uniqueness result for symmetric quasilinear hyperbolic systems (cf. \cite{tay}, section 16-17) that 
\begin{align} \label{upuq}
u_p = u_q \text{ in } \mathcal{U}_p \cap \mathcal{U}_q.
\end{align}
We now set 
\[M:= \cup_{p \in \Sigma} \mathcal{U}_p \subset \mathbb{R} \times \Sigma. \] 
By \eqref{upuq} the $g^{\mathcal{U}_p}$ define a global Lorentzian metric on $M$ on which $\partial_t$ is a timelike vector field by choice of the $\mathcal{U}_p$ in step 1. We equip $M$ with the time orientation induced by $\partial_t$. By the previous local constructions, $(\Sigma,g^{\Sigma})$ embeds into $(M,g)$ with Weingarten tensor $W$. Using $\partial_t$, the orthonormal frame bundle of $M$ is reduced to the structure group $SO(n) \subset SO^+(1,n)$ and via pullback the spin structure of $\Sigma$ naturally induces a spin structure on $M$. If restricted to $\mathcal{U}_p$, this spin structure coincides with the one considered in step 1. These observations allow us to define a global spinor field $\phi \in \Gamma(M,S_M)$ by demanding that
\begin{align}
\phi_{|\mathcal{U}_p} = \left[\widetilde{s_{a=0,...,n}[g^{\mathcal{U}_p}]},\widetilde{\phi}^{\mathcal{U}_p} \right],
\end{align}
where $s_{a=0,...,n}[g^{\mathcal{U}_p}]$ has been given in \eqref{rea} and  $\widetilde{}$  denotes the local lift to the spin structure of $M$. By \eqref{upuq} this is well-defined as all data coincide on the overlap $\mathcal{U}_p \cap \mathcal{U}_q$. By construction, $\phi$ restricts on $\Sigma$ to $\ph$ (with the identifications from the second section). Moreover, it has been shown in step 2 that $\nabla \phi_{|\mathcal{U}_p} = 0$ and $g(V,V)_{|\mathcal{U}_p} =  0$ for every $p \in \Sigma$, that is
\begin{align*}
\nabla \phi & \equiv 0, \\
g(V_{\phi},V_{\phi}) & = 0.
\end{align*}
From the second step also follows that $E=E^{g,h}=0$.
Finally, we show that $\Sigma \subset M$ is a spacelike Cauchy hypersurface. To this end, let $\gamma:I \rightarrow M = \cup_{p \in \Sigma} \mathcal{U}_p$ be an inextendible timelike curve and let $t^* \in I$ be any fixed parameter. Let $p \in \Sigma$ such that $\gamma(t^*) \in \mathcal{U}_p$. For such fixed $p$ we consider the restricted curve
\[ \gamma_{|\gamma^{-1}(\mathcal{U}_p)} : \gamma^{-1}(\mathcal{U}_p) \rightarrow \mathcal{U}_p, \]
which is an inextendible timelike curve in the globally hyperbolic manifold $(\mathcal{U}_p,g^{\mathcal{U}_p})$. Thus, the spacelike Cauchy hyersurface $\Sigma_p \subset \Sigma$ in $\mathcal{U}_p$ is met by $\gamma_{|\gamma^{-1}(\mathcal{U}_p)}$. It remains to show that $\gamma$ meets $\Sigma$ at most once. Wrt. the splitting $\mathbb{R} \times \Sigma$ we decompose
\[ \gamma = (\gamma_t, \gamma_{\Sigma}) \]
and compute
\begin{align*}
0 > g(\dot \gamma_t \partial_t, \dot \gamma_t \partial_t) + 2 \cdot g(\dot \gamma_t \partial_t, \dot \gamma_{\Sigma}) + g(\dot \gamma_{\Sigma}, \dot \gamma_{\Sigma}).
\end{align*}
Let us assume that there is $\tau \in I$ with $\dot \gamma_t(\tau) = 0$ Let $q:= \gamma(\tau)$. Then $0 > g^{\mathcal{U}_q}(\dot \gamma_{\Sigma}(\tau), \dot \gamma_{\Sigma}(\tau))$. This, however, contradicts the condition \eqref{3sa} imposed on $\mathcal{U}_q$ and $g^{\mathcal{U}_q}$ in step 1. Consequently, $\gamma_t:I \rightarrow \mathbb{R}$ is strictly monotone, and thus $\gamma_t = 0$ has at most one solution. In total $\gamma$ intersects $\Sigma$ exactly once. It follows that $(M,g)$ is globally hyperbolic with Cauchy hypersurface $\Sigma$ and parallel null spinor $\phi$. This finishes the proof of Theorem \ref{mother}.
$\hfill \Box$ \\
\begin{re}
It follows from the proof of Theorem \ref{mother} that if $\Sigma$ is actually compact, then $M$ can be chosen to be of the form $M = (-\epsilon,\epsilon) \times \Sigma$ for some $\epsilon$. Compact solutions to the constraint equations \eqref{spinor-1} and \eqref{spinor-2} are known to exist, cf. \cite{blli} and references therein.
\end{re}

%compare with results from analytic case

\section{Special cases}
First, we establish a connection between Theorem \ref{mother} and the results obtained in \cite{blli} for the Cauchy problem for parallel spinors in the analytic category. In the latter reference it was found that if the data $(\Sigma, g^{\Sigma},W,\ph)$ as appearing in Theorem \ref{mother} are all analytic, then one can for any prescribed analytic and positive lapse function $\lambda$ defined on $\mathbb{R} \times \Sigma$ find a unique family of analytic Riemannian metrics $g_t$ with $g_0 = g^{\Sigma}$ such that the analytic metric 
\begin{align}
g_{\text{analytic}} = -\lambda^2 dt^2 + g_t \label{stre}
\end{align}
 defined around $\Sigma$ admits a parallel isotropic spinor $\phi$ which restrict on $\Sigma$ to $\ph$. The proof is based on rewriting the conditions for $g_t$ and other data which follow from the existence of a parallel spinor in a system of Cauchy-Kowalewski type, which has no general solution theory in the smooth setting.\\
In order to make contact with Theorem \ref{mother}, let us as above assume that the initial data $(\Sigma, g^{\Sigma},W,\ph)$ are analytic. Then we use for some fixed analytic $\lambda$ the analytic metric \eqref{stre} from \cite{blli} as background metric in Theorem \ref{mother}, i.e. we set
\begin{align}
h = g_{\text{analytic}}. \label{sur} 
\end{align}
It follows immediately that for the resulting solution-metric $g_h$ as appearing in Theorem 1 we have $g_h=g_{\text{analytic}}$ because $g_{\text{analytic}}$ satisfies properties (1)-(3) from Theorem \ref{mother} by construction, property (4) is a trivial consequence of \eqref{sur} and these properties uniquely determine $g_h$. Thus, for analytic initial data, the metrics $g_h$ constructed via Theorem \ref{mother} include the metrics constructed in \cite{blli}.
However, for arbitrary smooth initial data it remains unclear whether one can specify a background metric $h$ in such a way that $g_h$ satisfies $\partial_t \perp T\Sigma$ as in \eqref{stre}.
%Given a - not necessarily globally hyperbolic - Ricci-flat Lorentzian manifold $(M,g)$ with spacelike Cauchy hypersurface $(\Sigma,\g)$ which embeds into $M$ with Weingartentensor $W$, the condition $\text{Ric}^{g} = 0$ imposes the following \textit{constraint conditions} on $(\Sigma,\g)$ (cf Prop...)

%\begin{align*}
 %\text{scal}^{\g} &=\text{tr}_{\g}(W^2) - \text{tr}^2_{\g} W, \\
  %\text{tr}_{\g}W  &= - \text{div}^{\g} W.
%\end{align*}
%The main result of this section is to show a special case of the converse direction, leading to a construction principle for Ricci-flat globally hyperbolic Lorentzian manifolds admitting a parallel spinor.

\begin{re}
\cite{blli} also solves the more general Cauchy problem for lightlike parallel vector fields in the analytic setting. We do not know whether this can be extended to the smooth category as well. The problem here is to find an analogue of \eqref{3}-\eqref{4} for the vector field case.
\end{re}
%However, we can not ensure that in the general smooth setting $g$ is always of type \eqref{analy}.

Returning to the smooth category, it is natural to ask under which additional conditions on the initial data $(\Sigma,g^{\Sigma},W,\ph)$ the manifolds $(M,g)$ constructed via Theorem \ref{mother} are actually Ricci-flat.
%restr spinorts lor spin manfiold ricci flat
 %Kor
\begin{mcor} \label{zur}
In the setting of Theorem \ref{mother}, the Lorentzian manifold $(M,g)$ is Ricci-flat if and only if the tensor $W$ additionally satisfies the constraint
\begin{align}
d \text{tr}_{g^{\Sigma}}W  + \delta^{\Sigma} W = 0. \label{co3}
\end{align}
\end{mcor}

\bprf
As a consequence of \eqref{leg}, condition \eqref{co3} is equivalent to $f_{|\Sigma} = 0$. In this case, $f$ can be set to zero in the system \eqref{system} which remains symmetric quasilinear hyperbolic. One then proceeds as in the proof of Theorem \ref{mother} with $f=0$, which is equivalent to $\text{Ric} = 0$.
\eprf

\begin{re}
The equations \eqref{co3} as well as \eqref{fsigma2} set to zero are precisely the constraint equations for the Cauchy problem for the vacuum Einstein equations.
\end{re}

\begin{re}
Symmetric $(1,1)$-tensors on $\Sigma$ satisfying $\eqref{co3}$ can be interpreted as \textit{generalized Codazzi tensors}. By definition, a $(1,1)$ tensor $W$ on $\Sigma$ is Codazzi if $(\nabla^{\Sigma}_X W)(Y)$ is symmetric in $X,Y \in T\Sigma$. Indeed, for $W$ Codazzi we find for $X \in T\Sigma$ and an orthonormal basis $(s_1,...,s_n)$ which is parallel in a fixed point that at this point 
\begin{align*}
{\delta}^{\Sigma}W(X) & = - \sum_i g^{\Sigma}((\nabla^{\Sigma}_{s_i}W)(X),s_i)\\
& = - \sum_i  g^{\Sigma}((\nabla^{\Sigma}_{X}W)(s_i),s_i) \\
 & = - X(\text{tr}_{g^{\Sigma}}W),
\end{align*}
i.e. $W$ satisfies \eqref{co3}. \cite{bam} shows that any complete Riemannian manifold $\Sigma$ carrying an imaginary generalized Killing spinor with $W$ being a uniformly bounded Codazzi tensor can be extended to a globally hyperbolic and Ricci flat Lorentzian manifold $\mathbb{R} \times \Sigma$ with parallel spinor. Thus, Corollary \ref{zur} is a generalization thereof to a class of more general imaginary $W-$Killing spinors. Note, however, that in contrast to \cite{bam} we cannot give an explicit description of $g$ on $M$.
\end{re}

\begin{bsp}
Any warped product $\Sigma:=  \mathbb{R} \times_h \cF$, where
$(\F,\g_{\F})$ is a $(n-1)$-dimensional complete Riemannian manifold, $h$ is a smooth positive function on $\mathbb{R}$ and $g^{\Sigma} = ds^2 + h^2 \cdot \g_{\F}$, admits the closed
conformal vector field $U(s,x) = h(s) \partial_s(s,x)$ of length
$u=h$, and the endomorphism $\W:= b \,\Id_{T\M}$ with $b:= -
\ln(h)'$ satisfies $\nabla U = -u\W$. If $(\cF,\g_{\cF})$ is spin
and has a parallel spinor field, then $\M= L\times_h \cF$ is spin as
well with an imaginary $\W= b\,\Id_{T\M}$--Killing spinor $\widetilde{\ph}$. For a
proof of this, see \cite{rad}.\\
It is easy to verify that $\W$ solves \eqref{co3} iff $b  =$\text{const}. However, we will show that for sufficiently generic $h$ Corollary \ref{zur} still applies to this situation by modifying the generalized Killing spinor.\\
To this end, we set $\ph := f \cdot \widetilde{\ph}$ for a function $f: \mathbb{R} \subset \mathbb{R} \times \cF \rightarrow \mathbb{R}$ to be determined. It is straightforward to compute using $\partial_s \cdot \widetilde{\ph} = i \cdot \widetilde{\ph}$ that
\begin{align}
 \nabla^{\Sigma}_X \ph = \frac{i}{2} \cdot \left(b \cdot X - 2 X(\ln(f)) \cdot \partial_s \right) \cdot \ph,
\end{align}
i.e. $\ph$ is an imaginary generalized $W_{f}$-Killing spinor with 
\begin{align}
W_{f} = b \cdot \text{Id}_{T\Sigma} - 2 d(\ln(f)) \otimes \partial_s. \label{wf}
\end{align}
Note that $W_{f}$ is symmetric due to the special form of $g^{\Sigma}$. Moreover, $\ph$ satisfies the algebraic constraint \eqref{spinor-2} as $\widetilde{\ph}$ does. $d \text{tr}_{g^{\Sigma}}W_{f}(X)  + \delta^{\Sigma} W_{f}(X)$ holds automatically for all $X$ tangent to $\cF$ and inserting $X=\partial_s$ leads to the condition
\begin{align}
b \cdot \ln(f)' \stackrel{!}{=} \frac{1}{2n}(n-1) b'. \label{lel}
\end{align}
Thus, if $h$ is chosen such that $f:= (-\ln(h)')^{\frac{n-1}{2n}}$ is defined on all of $\mathbb{R}$, \eqref{lel} holds and by Corollary \ref{zur} $(\F,\g_{\F})$ can be embedded in a Ricci-fat Lorentzian manifold admitting a parallel spinor $\phi$ (arising from extending $\ph$ not $\widetilde{\ph}$!). Note that for generic choice of $h$, $W_{f}$ is not Codazzi. Moreover, for $h = s^{\frac{2(n-1)}{2n}}$ (and defining the warped product on $I \times \cF$ for appropriate $I \subset \mathbb{R}$), \eqref{wf} implies $W_{f}(\partial_s) = 0$, i.e. $W_{f}$ is not invertible.
\end{bsp}

%examples of ricci flat which are not codazzo

\textbf{Acknowledgements.} The author would like to thank Helga Baum and Thomas Leistner for
helpful discussions and useful comments and on a first draft of the paper.

\small
\bibliographystyle{plain}
\bibliography{literatur}

%\bibliographystyle{abbrv}
%%\bibliographystyle{amsalpha}
%\bibliography{GEOBIB}
%\end{document}
\end{document}